\documentclass[sn-mathphys,Numbered]{A_maximal_oscillatory_operator_on_compact_manifolds}

\usepackage{graphicx}%
\usepackage{multirow}%
\usepackage{amsmath,amssymb,amsfonts}%
\usepackage{amsthm}%
\usepackage{mathrsfs}%
\usepackage[title]{appendix}%
\usepackage{xcolor}%
\usepackage{textcomp}%
\usepackage{manyfoot}%
\usepackage{booktabs}%
\usepackage{algorithm}%
\usepackage{algorithmicx}%
\usepackage{algpseudocode}%
\usepackage{listings}%

\theoremstyle{thmstyleone}%
\theoremstyle{thmstyletwo}%

\theoremstyle{thmstylethree}%

\begin{document}

\title[Article Title]{ A maximal oscillatory operator on compact manifolds}


\author*[1]{\fnm{Ziyao} \sur{Liu}}\email{zy.liu@zjnu.edu.cn}

\author[1]{\fnm{Jiecheng} \sur{Chen}}\email{jcchen@zjnu.edu.cn}

\author[1,2]{\fnm{Dashan} \sur{Fan}}\email{fandashan2@zjnu.edu.cn}

\affil*[1]{\orgdiv{Department of Mathematical Science}, \orgname{Zhejiang Normal University}, \orgaddress{\city{Jinhua}, \postcode{321004}, \state{Zhejiang}, \country{CHINA}}}

\affil[2]{\orgdiv{Department of Mathematical Science}, \orgname{University of Wisconsin-Milwaukee}, \orgaddress{\city{Milwaukee}, \postcode{WI 53201}, \state{Wisconsin}, \country{USA}}}



\abstract{This is a continuation of our previous research about an oscillatory integral operator $T_{\alpha, \beta}$ on compact manifolds $\mathbb{M}$. We prove the sharp $H^{p}$-$L^{p,\infty}$ boundedness on the maximal operator $T^{*}_{\alpha, \beta}$ for all $0<p<1$. As applications, we first prove the sharp $H^{p}$-$L^{p,\infty}$ boundedness on the maximal operator corresponding to the Riesz means $I_{k,\alpha}(|\mathcal{L}|)$ associated with the Schrödinger type group $e^{is\mathcal{L}^{\alpha/2}}$ and obtain the almost everywhere convergence of $I_{k,\alpha}(|\mathcal{L}|)f(x,t)\to f(x)$ for all $f\in H^{p}$. Also, we are able to obtain the convergence speed of a combination operator from the solutions of the Cauchy problem of fractional Schrödinger equations. All results are even new on the n-torus $T^{n}$.}



\maketitle
\section{Introduction}\label{sec1}

Let $n\geq2$, and $(\mathbb{M}, g)$ be an $n$-dimensional compact manifold without boundary. In this article, we study an oscillating operator
$T_{\alpha, \beta, t}$ on $\mathbb{M}$ and obtain its sharp maximal estimate, where $0<\alpha<1$ and $\beta>0$.

We now start introducing the definition of operator $T_{\alpha, \beta, t}$. The Laplacian $\mathcal{L}$ on $\mathbb{M}$ is given locally by
$$\mathcal{L}=\frac{1}{\sqrt{\det\left[g_{i, j}\right]}} \sum_{i, j=1}^n \frac{\partial}{\partial x_i}\left(g^{i, j} \sqrt{\det\left[g_{i, j}\right]} \frac{\partial}{\partial x_j}\right),$$
where $g_{i, j}=g\left(\frac{\partial}{\partial x_i}, \frac{\partial}{\partial x_j}\right)$, $\left[g_{i, j}\right]$ is a real and positive definite matrix and $\left[g^{i, j}\right]$ is the inverse matrix of $\left[g_{i, j}\right]$.
The full set of eigenvalues of $\mathcal{L}$ is given by
$$\left\{\lambda_1^2, \lambda_2^2, \lambda_3^2, \ldots\right\},$$
where  $\lambda_1^2 \leq \lambda_2^2 \leq \lambda_3^2 \leq \ldots \rightarrow \infty$; and  the family of their corresponding eigenfunctions $\left\{e_{j}(x)\right\}$ forms a real orthonormal basis of $L^{2}(\mathbb{M})$ (see \cite{bib9,bib10,bib20}). Without loss of generality, we assume $e_{1}\left( x\right) =1$ if $\lambda
_{1}=0$. For any $f\in{L^{2}(\mathbb{M})}$, its Fourier expansion is of the form
$$
f(x)=\sum_{j}\left\langle f, e_j\right\rangle e_j(x)=\sum_{j} e_j(f)(x),
$$
where $$e_j(f)(x)=\left\langle f, e_j\right\rangle e_j(x)$$
and $$\left\langle f, e_j\right\rangle=\int_{\mathbb{M}}e_{j}(y)f(y)dy.$$

Let $\Phi$ be a $C^{\infty}$ nonnegative function on the real line $\mathbb{R}$ that satisfies $\Phi(\lambda)\equiv1$ if $\lvert\lambda\rvert\ge2$ and $\Phi(\lambda)\equiv0$ if $\lvert\lambda\rvert\le1.$ For fixed positive numbers $\alpha$ and $\beta$, the oscillating multiplier operator that we mainly study is defined as 
$$
T_{\alpha, \beta, t}(f)(x)=\sum_{j} \frac{e^{i\lvert t \lambda_j\rvert^\alpha}}{\lvert t \lambda_j\rvert^\beta} \Phi\left(\lvert t \lambda_j\rvert\right)\left\langle f, e_j\right\rangle e_j(x),
$$
where  the parameter $t$ is positive. In the above formulation, we may initially assume $f\in C^{\infty}(\mathbb{M})$. Thus, $T_{\alpha, \beta, t}(f) $ can be written in its integral form as

$$T_{\alpha, \beta, t}(f)(x)=\int_{\mathbb{M}} \Omega_{\alpha, \beta}(x, y, t) f(y)dy,
$$
where 
$$
\Omega_{\alpha, \beta}(x, y, t)=\sum_{j} \frac{e^{i\lvert t \lambda_j\rvert^\alpha}}{\lvert t \lambda_j\rvert^\beta} \Phi\left(\lvert t \lambda_j\rvert\right)e_j(x)e_j(y)
$$
is the distributional kernel.

In order to obtain the almost everywhere convergence theorem, we study the maximal operator
$$
T_{\alpha, \beta}^*(f)(x)=\sup _{0<t \leq \sigma}\left|T_{\alpha, \beta, t}(f)(x)\right|=\sup _{0<t \leq \sigma}\left|\int_{\mathbb{M}} \Omega_{\alpha, \beta}(x, y, t) f(y) d y\right|,
$$
where $\sigma$ is any fixed positive number satisfying $\sigma<r_{0}/6$, where $r_{0}$ is the injectivity radius of $\mathbb{M}$. To illustrate $r_{0}$, we let $Exp_{z}$ be the exponential map at $z\in\mathbb{M}$. Due to the compactness of $\mathbb{M}$, there exists an $r_{z}>0$ such that for every $z\in\mathbb{M}$, $Exp_{z}$ is injective from the ball $B(0, r_{z})$ in the tangent space $T_{z}$ at $z$ into the ball $B(z, r_{z})$ in $\mathbb{M}$. The largest $r_{z}$, denote it $r_0$, with this property is called the injectivity radius of $\mathbb{M}$ (see \cite{bib6}).

It is worth to remark that $T_{\alpha, \beta}^*$ is a sort of pseudo-differential operator
$$T_{\alpha, \beta}^*(f)(x)=\sum_{j}\frac{e^{i\lvert t(x) \lambda_j\rvert^\alpha}}{\lvert t(x) \lambda_j\rvert^\beta} \Phi\left(\lvert t(x) \lambda_j\rvert\right)\left\langle f, e_j\right\rangle e_j(x),$$
where $t(x)$ is a function of $x$ and takes values in $(0, \sigma] $.

Let $L^{p}(\mathbb{M})$ $(1\le p\le\infty)$ be Lebegue spaces and $L^{p,\infty}(\mathbb{M})$ be weak $L^{p}(\mathbb{M})$ spaces. It is well known that the Hardy space $H^{p}(\mathbb{M})$ is a good substitute of $L^{p}(\mathbb{M})$ when $0<p\le1$ and $L^{p}(\mathbb{M})=H^{p}(\mathbb{M})$ for $1<p<\infty$. We recall that the Hardy space $H^{p}(\mathbb{M})$ is  the set of all distributions $f$ satisfying (see \cite{bib10})
$$
\|f\|_{H^p(\mathbb{M})}=\left\|\sup _{0<t<\infty}\left|e^{-t\mathcal{L}} f\right|\right\|_{L^p(\mathbb{M})}<\infty,
$$
where $e^{-t\mathcal{L}}$ is the heat operator defined on distributions $f$ by
$$e^{-t\mathcal{L}}(f)(x)=\sum_{j}e^{-t\lambda_{j}^{2}}\left\langle f, e_j\right\rangle e_j(x).$$

We also need to use the atomic characterization of $H^{p}(\mathbb{M})$ $(0<p<1)$. For $0<p\le 1$, any $f\in H^{p}(\mathbb{M})$ admits an atomic decomposition (see \cite{bib9,bib10})
$$f=\sum_{k}c_{k}a_{k},$$
where $\left\{c_{k}\right\}$ is a sequence of complex numbers satisfying (For the definition of $\preceq$, see page 7.)
$$\sum_{k}\lvert c_{k}\rvert^{p}\preceq\left\Vert f \right\Vert^{p}_{H^{p}(\mathbb{M})}, $$
and these functions $a_{k}$ are either regular or exceptional atoms. From the previous definition of $r_{0}$, we may assume that any ball in $\mathbb{M}$ have radius smaller than $r_{0}/10$. Given a ball $B(z,r)\subset\mathbb{M}$, we can define  polynomials on $B(z,r)$ as the images of polynomials on the tangent space $T_{z}$ via the exponential map. The regular $p$-atom $a$ (or $p$-atom) is a function in $L^{2}(\mathbb{M})$ satisfying the following conditions:  

$(i)$ the support of $a$ is contained in a ball $B(z, r)$, 

$(ii)$ $\left\Vert a\right\Vert_{L^{2}(\mathbb{M})}\le\lvert B(z, r)\rvert^{1/2-1/p}$, where $\lvert B\rvert$ denotes the measure of a set $B\subset\mathbb{M}$,

$(iii)$ $\int_{\mathbb{M}}a(x)\mathcal{Q}(x)dx=0$,
for every polynomial $\mathcal{Q}$ on $B(z, r)$ of degree at most $\left[n(1/p-1)\right]$. Condition $(iii)$ is called the cancellation condition of $a$.

In addition, an exceptional atom $a$ is an $L^{\infty}$ function on $\mathbb{M}$ satisfying
$$\left\Vert a\right\Vert_{L^{\infty}(\mathbb{M})}\le1.$$

The Riesz potential $\mathfrak{R}_{s}$ is defined via the Fourier expansion by
$$
\mathfrak{R}_s(f)(x)=\sum_{j}\left|\lambda_j\right|^s e_j(f)(x)
$$
when $s>0$ and 
$$
\mathfrak{R}_s(f)(x)=\sum_{j \neq 1}\left|\lambda_j\right|^s e_j(f)(x)
$$
when $s<0$.  The Sobolev imbedding theorem states that (see \cite{bib9,bib11,bib13})
$$\left\Vert \mathfrak{R}_{-s}(f) \right\Vert_{H^{p}(\mathbb{M})}\preceq \left\Vert f \right\Vert_{H^{q}(\mathbb{M})}$$
for any $s\ge0$, and $0<p, q<\infty$ satisfying
$$1/p=1/q-s/n.$$

Also, the space of all distributions $f$ satisfying
$$\left\|f \right\|_{H^{p}_{s}(\mathbb{M})}= \left\Vert \mathfrak{R}_{s}(f) \right\Vert_{H^{p}(\mathbb{M})}<\infty$$
is called the Sobolev $H^p$ space of order $s$, $s>0$.

The main purpose of this article is to establish the boundedness of $T_{\alpha, \beta}^*$ from $H^{p}(\mathbb{M})$ to $L^{p,\infty}(\mathbb{M})$ for $0<p<1$ under the optimal relation among $\alpha$, $\beta$ and $p$. This research topic has a deep background rooted in Euclidean space.  Precisely, the formulation of $T_{\alpha, \beta, t}$  is an analogue of the oscillating multiplier operator $S_{\alpha, \beta, t}$ on $\mathbb{R}^{n}$, which is defined as
$$S_{\alpha, \beta, t}(f)(x)=\int_{\mathbb{R}^{n}}\widehat{f}(\xi)\Phi (\lvert t\xi \rvert) \frac{e^{i\lvert t \xi\rvert^\alpha}}{\lvert t \xi\rvert^\beta} e^{i<x, \xi>}d\xi,$$
where $\Phi$ is the same as in the definition of $T_{\alpha, \beta, t}$ and $\widehat{f}$ denotes the Fourier transform of $f$ given by
$$\widehat{f}(\xi)=\int_{\mathbb{R}^{n}}f(x) e^{i<x, \xi>}d\xi.$$
When $t=1$, we can simply write it as
$$S_{\alpha, \beta}(f)(x)=S_{\alpha, \beta, 1}(f)(x).$$

    The study of the operator $S_{\alpha, \beta}$ in harmonic analysis is of interest and, on the other hand, it plays an important role in the study of Cauchy problems of Schrödinger and wave equations (see \cite{bib12,bib24}). We recall the following result obtained by Wainger \cite{bib24}, Fefferman and Stein \cite{bib12}, Sjölin \cite{bib21} and Miyachi \cite{bib18} in the classical case.

\textbf{Theorem A.} \textit{Assume} $\alpha\neq1$ \textit{and} $0<p<\infty$. $S_{\alpha, \beta}$ \textit{is bounded from} $H^{p}(\mathbb{R}^{n})$\textit{ to} $L^{p}(\mathbb{R}^{n})$ \textit{if and only if} 
$$|1/2-1/p| \le (n\alpha)^{-1}\beta.$$

On the other hand, the study of $S_{\alpha, \beta}$ on $\mathbb{R}^{n}$ was extended and generalized to many underlying spaces of different settings, see \cite{bib1,bib3,bib4,bib14,bib17}, among many others. Particularly, we recall the following theorem related to this research. 

\textbf{Theorem B.} \cite{bib6} \textit{Let} $n\ge2$ \textit{and let }$(\mathbb{M}, g)$ \textit{be an n-dimensional compact connected manifold without boundary. Assume} $0 <\alpha< 1$. \textit{Then} 

(i) $T^{*}_{\alpha, \beta}$\textit{ is bounded from } $L^{1}(\mathbb{M})$\textit{ to} $L^{1, \infty}(\mathbb{M})$ if  
$\beta>(n\alpha)/2.$

(ii) $T^{*}_{\alpha, \beta}$ \textit{is bounded on }$L^{p}(\mathbb{M})$ \textit{if} $$\max\left\lbrace\dfrac{2n\alpha-2\alpha}{2\beta+n\alpha-2\alpha}, 1 \right\rbrace <p<\dfrac{2n\alpha-2\alpha}{\min\left\lbrace0, n\alpha-\beta \right\rbrace }.$$

\textit{Additionally}, $T^{*}_{\alpha, \beta}$ \textit{is bounded on} $L^{2}(\mathbb{M})$\textit{ if }$\beta>\alpha/2$. 

(i) in Theorem B is even new on $\mathbb{R}^{n}$ or torus $T^{n}$. In fact, the result is true on $T^{n}$, since $T^{n}$ is a special case of $\mathbb{M}$.  Then a transference theorem in \cite{LF} gives the same result on $\mathbb{R}^{n}$. Clearly, comparing Theorem B with Theorem A, we find that more regularity is needed to  guarantee the $L^{p}$ boundedness of $T^{*}_{\alpha, \beta}$. This is quite natural because that $T^{*}_{\alpha, \beta}$ is a larger operator than $T_{\alpha, \beta}$. In Theorem B, the $L^{p}$ boundedness of $T^{*}_{\alpha, \beta}$ when $1<p\le2$ needs the assumption $\beta>n\alpha (\frac{1}{p}-\frac{1}{2})+\alpha(1-\frac{1}{p})$. This means that  when $1<p\le2$ we need an extra reguarity $\alpha(1-\frac{1}{p})$ comparing to the non-maximal operator. Similarly, we need an extra $\frac{\alpha}{p}$ in the case $p>2$. What is the sharpness of $\beta$ still remains open when $p>1$, even on the torus or on $\mathbb{R}^{n}$. However, when $0<p<1$, we obtain the  following sharp result, which is even new on $T^{n}$ based on our best knowledge. 

\textbf{Theorem 1.} \textit{Let} $n\ge2$ \textit{and let} $(\mathbb{M}, g)$ \textit{be an n-dimensional compact connected manifold without boundary. Assume} $0 <\alpha< 1$ \textit{and} $0<p<1$. \textit{Then we have the following statements.} 

(i) $T^{*}_{\alpha, \beta}$ \textit{is bounded from} $H^{p}(\mathbb{M})$ \textit{to} $L^{p,\infty}(\mathbb{M})$ \textit{if and only if}  
$$\beta\ge n\alpha(1/p-1/2) .$$

(ii) $T^{*}_{\alpha, \beta}$ \textit{is bounded from} $H^{p}(\mathbb{M})$ \textit{to} $L^{p}(\mathbb{M})$ \textit{if}  
$$\beta> n\alpha(1/p-1/2) .$$

\textbf{Remark} In fact, our results are sharp in sense that the non-maximal operator 
$T_{\alpha, \beta}$ is bounded on $H^{p}(\mathbb{M})$ if and only if $\beta\geq n\alpha(1/p-1/2)$.

Combining Theorem 1 and Theorem B in \cite{bib6}, we give a complete solution to the $H^{p}(\mathbb{M})$ boundedness of $T^{*}_{\alpha, \beta}$ when $0<p\le1$. The complete solution to the case $p>1$ is interesting and challenging, and still remains open. We hope to return this topic in our future research.

Theorem 1 can be applied in the study of Approximation Theory. We consider the Cauchy problem of the fractional Schrödinger equation, for any natural number $k$,
\begin{equation}
	\left\{
	\begin{array}{ll}
		\frac{\partial}{\partial t}u_{k}(x, t)-ik\mathcal{L}^{\alpha/2}u_{k}(x, t)= 0\\
		u_{k}(x, 0)=f(x).
	\end{array}
	\right.\notag
\end{equation}
The formal solution is given by
$$u_{k}(x, t)=A_{kt}^{\alpha}(f)(x)=\sum_{j}e^{ikt|\lambda_{j}|^{\alpha}}\left\langle f, e_j\right\rangle e_j(x).$$
It is an interesting problem to study the convergence speed of the combination operator
$$\lim_{t\to0}\sum_{k=1}^{N}c_{k}A_{kt}^{\alpha}(f)(x)=f(x), a.e.$$
related to the regularity of $f$ for an appropriate sequence $(c_{1}, c_{2},\dots, c_{N})$.

Now we have the following theorem.

\textbf{Theorem 2.} \textit{Fix} $0<p<1$ \textit{and} $0<\alpha<1$. \textit{If} $f\in H_{\beta}^{p}(\mathbb{M})$ with $\beta\ge n\alpha(1/p-1/2)$, \textit{then for any} $N>[\beta/\alpha]$ \textit{there exists a sequence of constants} $\{c_{k}\}_{k=1}^{N}$, \textit{such that}
$$\sum_{k=1}^{N}c_{k}A_{kt}^{\alpha}(f)(x)-f(x)=o(t^{\frac{\beta}{\alpha}})\quad \textup{as}\ t\to0$$ 
\textit{for almost all} $x\in \mathbb{M}.$

Here, $[\beta/\alpha]$ denotes the largest integer part of $\beta/\alpha$. We may obtain a similar result by applying Theorem 1 to the maximal operator corresponding to the Riesz means associated with the Schrödinger type group $e^{is\mathcal{L}^{\alpha/2}}$.

Let
$$I_{k,\alpha}(\mathcal{L})(f)(x,t)=kt^{-k}\int_{0}^{t}(t-s)^{k-1}e^{is\mathcal{L}^{\alpha/2}}(f)ds,\quad 0<\alpha<1,\ k>0,$$
and
$$I^{*}_{k,\alpha}(\mathcal{L})(f)(x)=\sup_{0<t\le\sigma}\left| kt^{-k}\int_{0}^{t}(t-s)^{k-1}e^{is\mathcal{L}^{\alpha/2}}(f)ds\right|,$$
where 
$$e^{is\mathcal{L}^{\alpha/2}}(f)(x)=\sum_{j}e^{is|\lambda_{j}|^{\alpha}}\left\langle f, e_j\right\rangle e_j(x).$$
We have the following theorem.

\textbf{Theorem 3.} \textit{Let} $n\ge2$ \textit{and let} $(\mathbb{M}, g)$ \textit{be an n-dimensional compact connected manifold without boundary. Assume} $0 <\alpha< 1$ \textit{and} $0<p<1$.  \textit{Then} $I^{*}_{k, \alpha}(\mathcal{L})$ \textit{is bounded from} $H^{p}(\mathbb{M})$ \textit{to} $L^{p}(\mathbb{M})$ \textit{if}  
$$k> n\alpha(1/p-1/2) .$$
\textit{Also}, $I^{*}_{k, \alpha}(\mathcal{L})$ \textit{is bounded from} $H^{p}(\mathbb{M})$ \textit{to} $L^{p, \infty}(\mathbb{M})$ \textit{if}  
$$k\geq n\alpha(1/p-1/2) .$$
\textit{In addition},
$$\lim_{t\to0}I_{k,\alpha}(\mathcal{L})(f)(x,t)=f(x)$$
\textit{almost everywhere for any} $f\in H^{p}(\mathbb{M})$, \textit{provided} $k\geq n\alpha(1/p-1/2) .$

This paper is organized as follows. In the second section, we represent some preliminary knowledge on $\mathbb{M}$. We will decompose $T^{*}_{\alpha, \beta}$ as a sum of two maximal operators $T^{*}_{\alpha, \beta, loc}$ and $T^{*}_{\alpha, \beta, \infty}$. The estimate on the kernel of $T^{*}_{\alpha, \beta, loc}$ will be obtained in the third section. In Section 4, we shall use the Sobolev imbedding theorem to show the $H^{p}$ boundedness of $T^{*}_{\alpha, \beta, \infty}$. In Section 5, we prove Theorem 1. The proofs of Theorem 2 and Theorem 3 can be found in Section 6.

Throughout the article, we use the notation $A \preceq B$ to mean that there is a positive constant $C$ independent of all essential variables such that $A \le CB$ ($\succeq$ similarly). The notation $A \simeq B$ denotes $A \preceq B$ and $A \succeq B$, 
and the notation $A \approx B$ means that there exists a constant C such that $A=CB$. We also use the symbol [$r$] to denote the largest integer part of a real number $r$.

\section{Preliminary Knowledge and Some Known Results}\label{sec2}

Throughout this paper, we denote the Riemannian distance between two points $%
x,y\in \mathbb{M}$ by $|x-y|$. Write 
\begin{equation*}
	\mu _{\alpha ,\beta }(\lambda ,t)=\frac{e^{i\lvert t\lambda \rvert ^{\alpha
	}}}{\lvert t\lambda \rvert ^{\beta }}\Phi (\lvert t\lambda \rvert ),
\end{equation*}
$\widehat{\mu }_{\alpha ,\beta }$ is its Fourier cosine transform on the $%
\lambda $ variable (see Section 3.1) and $\mu _{\alpha ,\beta }(\lambda ,1)=\mu _{\alpha
	,\beta }(\lambda )$. By the spectral theory (see \cite{bib2,bib16}), the kernel $\Omega
_{\alpha ,\beta }(x,y,t)$ can be written in the expression 
\begin{equation*}
	\begin{aligned} 
	\Omega_{\alpha, \beta}(x, y, t) =&\frac{2}{\pi
			t}\sum_{j}\left(\int_{0}^{\infty} \widehat {\mu}_{\alpha, \beta} (\tau
		/t) \cos(|\lambda_{j} |\tau)d\tau \right)e_{j}(x)e_{j}(y)\\ =&\frac{2}{\pi t}
		\int_{0}^{\infty} \widehat {\mu}_{\alpha, \beta}(\tau/t)\cos (\tau
		\sqrt{\mathcal{L}})(x, y)d\tau, 
		\end{aligned}
\end{equation*}%
where $\cos (\tau \sqrt{\mathcal{L}})(x,y)$ is the distribution kernel
associated with the operator $\cos (\tau \sqrt{\mathcal{L}})$ defined by 
\begin{equation*}
	\cos (\tau \sqrt{\mathcal{L}})f=\sum_{j}\cos (|\lambda _{j}|\tau
	)e_{j}(f).
\end{equation*}%
By the Hadamard construction of a parametrix (see \cite{bib15} or
\cite{bib1,bib2,bib17}), we have a $\delta >0$ such that if $|\tau |<\delta $ then 
\begin{equation}
	\cos (\tau \sqrt{\mathcal{L}})(x,y)=\sum_{j=0}^{\infty }U_{j}(x,y)|\tau |%
	\frac{\left( \tau ^{2}-|x-y|^{2}\right) _{+}^{j-\frac{n+1}{2}}}{\Gamma
		\left( j-\frac{n-1}{2}\right) },  \label{1}
\end{equation}%
where $U_{0}(x,y)\neq 0$ and the functions $U_{j}(x,y)\in C^{\infty }(%
\mathbb{M}\times \mathbb{M})$ satisfy
\begin{equation}
	\sum_{j=0}^{\infty }|\frac{\partial ^{\gamma }}{\partial y^{\gamma }}%
	U_{j}(x,y)|<C,  \label{2}
\end{equation}%
uniformly for $(x,y)\in \mathbb{M}\times B(x,\delta )$ and all multi-indices 
$\gamma $. With a little abuse we use the notation $\frac{\partial^\gamma}{\partial y^\gamma}$ to denote the $\gamma$-th derivative at the point $y$ in $\mathbb{M}$ in some local coordinates (see \cite{bib9}).

In the above asymptotic expansion in \eqref{1}, we observe the finite
propagation speed property 
\begin{equation*}
	\cos (\tau \sqrt{\mathcal{L}})(x,y)=0\quad \textup{if}\ |x-y|\geq |\tau |.
\end{equation*}%
Also, from \eqref{1} we agree that when $s\geq 1$ the distribution $\frac{%
	\tau _{+}^{-s}}{\Gamma (-s+1)}$ is defined recursively by 
\begin{equation}
	\int_{0}^{\infty }\frac{\tau _{+}^{-s}}{\Gamma (-s+1)}f(\tau )d\tau
	=-\int_{0}^{\infty }\frac{\tau _{+}^{-s+1}}{\Gamma (-s+2)}\frac{d}{d\tau }%
	f(\tau )d\tau . \label{3}
\end{equation}%
Without loss of generality, we assume $r_{0}\leq \delta $. Then, in the
asymptotic expansion of $\cos (\tau \sqrt{\mathcal{L}})(x,y)$, we may assume
that the inequality in \eqref{2} holds uniformly on $(x,y)\in \mathbb{M}%
\times \mathbb{M}$.

For the injectivity radius $r_{0}$ of $\mathbb{M}$, we fix an even $C^{\infty}$ function $\psi$ which satisfies $0\le\psi(s)\le1$, and $\psi(s)=1$ if $|s|\le r_{0}/3$ and $\psi(s)=0$ if $|s|\ge 2r_{0}/3$. Set
$$\varphi(s)=1-\psi(s).$$
We decompose 
$$\Omega_{\alpha, \beta}(x, y, t)=\Omega_{\alpha, \beta, loc}(x, y, t)+\Omega_{\alpha, \beta, \infty}(x, y, t),$$
where 
$$\Omega_{\alpha, \beta, loc}(x, y, t)=\frac{1}{\pi t} \int_{0}^{\infty} \widehat {\mu}_{\alpha, \beta}(\tau/t)\psi(\tau)\cos (\tau \sqrt{\mathcal{L}})(x, y)d\tau$$ 
and 
$$\Omega_{\alpha, \beta, \infty}(x, y, t)=\frac{1}{\pi t} \int_{0}^{\infty} \widehat {\mu}_{\alpha, \beta}(\tau/t)\varphi(\tau)\cos (\tau \sqrt{\mathcal{L}})(x, y)d\tau.$$
The corresponding operators, respectively, are
$$T_{\alpha, \beta, loc}(f)(x, t)=\int_{\mathbb{M}}\Omega_{\alpha, \beta, loc}(x, y, t)f(y)dy,$$ 
$$T_{\alpha, \beta, \infty}(f)(x, t)=\int_{\mathbb{M}}\Omega_{\alpha, \beta, \infty}(x, y, t)f(y)dy,$$
and the corresponding maximal operators, respectively, are defined as
$$T^{*}_{\alpha, \beta, loc}(f)(x)=\sup_{0<t\le\sigma}|T_{\alpha, \beta, loc}(f)(x, t)| $$
and
$$T^{*}_{\alpha, \beta, \infty}(f)(x)=\sup_{0<t\le\sigma}|T_{\alpha, \beta, \infty}(f)(x, t)| .$$

\section{Estimate on Kernel of $T_{\alpha, \beta, loc}$}\label{sec3}

\subsection{Derivative estimate of $\widehat{\mu}_{\alpha, \beta}(\tau)$}\label{subsec2}

In order to obtain a precise derivative estimate of $\Omega _{\alpha ,\beta
	,loc}(x,y,t)$, in this subsection we aim to study the derivative estimate of 
$\widehat{\mu }_{\alpha ,\beta }(\tau )$, which is the Fourier cosine
transform of $\mu _{\alpha ,\beta }$, 
$$\widehat{\mu}_{\alpha, \beta}(\tau)=2\int_{0}^{\infty}\frac{e^{i \lambda^\alpha}}{ \lambda^\beta}\Phi ( \lambda)\cos(\tau\lambda)d\lambda,$$
where $\Phi$ is the function in the definition of $T_{\alpha, \beta, t}$. 

Let $\phi$ be a $C^{\infty}$ nonnegative even function supported in the interval $[1/2, 2]\cup[-2, -1/2]$. Pick another $C^{\infty}(\mathbb{R})$ even function $\Psi_{0}$ satisfying that $\Psi_{0}(\lambda)\equiv1$ on the set $\left\{\lambda: |\lambda|\le1/2\right\}$ and $\Psi_{0}$ is supported on the set $\left\{\lambda: |\lambda|\le1\right\}$. Also, $\phi$ and $\Psi_{0}$ satisfying 
$$\sum_{k=0}^{\infty}\phi(\frac{|u|}{2^{k}})+\Psi_{0}(|u|)\equiv1$$
for any $u\in\mathbb{R}$.

Without loss of generality, we may write
$$\mu_{\alpha, \beta}(\lambda)=\sum_{k=0}^{\infty}\mu_{\alpha, \beta, k}(\lambda),$$
where 
$$\mu_{\alpha, \beta, k}(\lambda)=\frac{e^{i\lvert  \lambda\rvert^\alpha}}{\lvert \lambda\rvert^\beta}\phi (\frac{\lvert \lambda \rvert}{2^{k}}).$$
In the above formula, by the definition we may assume $\Psi_{0}=0$ and $\Phi\equiv1$ when they appear in the formula of $\mu_{\alpha, \beta, k}$ without loss of generality.

By this definition, we have that
$$\widehat{\mu}_{\alpha, \beta}(\tau)=\sum_{k=0}^{\infty}\widehat{\mu}_{\alpha, \beta, k}(\tau),$$
where 
$$\widehat{\mu}_{\alpha, \beta, k}(\tau)=2\int_{0}^{\infty}\frac{e^{i \lambda^\alpha}}{ \lambda^\beta}\phi (\frac{ \lambda}{2^{k}})\cos(\tau\lambda)d\lambda.$$

Choose $C^{\infty}$ nonnegative functions $\Gamma_{m}, m=1, 2, 3,$
satisfying
$$\textup{supp}\ \Gamma_{1}\subset\left\{\tau\in\mathbb{R}:|\tau|\le c_{1}\right\}\  \textup{and}\ \Gamma_{1}\equiv 1\ \textup{if}\ |\tau|\le c_{1}/2,$$
$$\textup{supp}\ \Gamma_{2}\subset\left\{\tau\in\mathbb{R}:|\tau|\ge c_{2}\right\}\  \textup{and}\ \Gamma_{2}\equiv 1\ \textup{if}\ |\tau|\ge 2c_{2},$$
and
$$\Gamma_{3}(\tau)=1-\Gamma_{1}(\tau)-\Gamma_{2}(\tau),$$
where we fix $c_{2}$ larger than $n2^{100}\alpha^{-1}$ and let $c_{1}=c_{2}^{-1}.$

In addition, for $m=1, 2, 3$ and $k=1, 2, \dots,$ we write
$$\Gamma_{m, k}(\tau)=\Gamma_{m}(2^{k(1-\alpha)}\tau).$$
Note that $\Gamma_{m, k},\ m=1, 2, 3,$ are supported, respectively, in the sets 
$$E_{1, k}=\{\tau:|\tau|\le c_{1}2^{k(\alpha-1)}\},$$
$$E_{2, k}=\{\tau:|\tau|\ge c_{2}2^{k(\alpha-1)}\},$$
$$E_{3, k}=\{\tau:c_{1}2^{k(\alpha-1)-1}\le |\tau|\le c_{2}2^{k(\alpha-1)+1}\}.$$
By the partition of the unity, we obtain that
$$\frac{\partial^{L}}{\partial \tau^{L}}\widehat{\mu}_{\alpha, \beta, k}(\tau)=\sum_{m=1}^{3}\Gamma_{m, k}(\tau)(\frac{\partial^{L}}{\partial \tau^{L}}\widehat{\mu}_{\alpha, \beta, k}(\tau)).
$$

Hence, we have the following estimates. Since the proof of the following Lemma 3.1 is similar to the proof of Lemma 1 in paper \cite{bib7}, we omit it.

\bigskip
\textbf{Lemma 3.1.} \textit{Let} $0<\alpha<1$ \textit{and} $\beta>0$. \textit{For any nonnegative integer} $L$, \textit{we have that, if} $|\tau|\le 200$, 
$$
\left|\frac{\partial^L}{\partial \tau^L} \widehat{\mu_{\alpha, \beta}}(\tau)\right| \preceq \max \left\{|\tau|^{\frac{-L}{1-\alpha}+\frac{\alpha-2+2 \beta}{2(1-\alpha)}}, 1\right\},
$$
\textit{and if} $|\tau|\ge 200$,
$$
\left|\frac{\partial^L}{\partial \tau^L} \widehat{\mu_{\alpha, \beta}}(\tau)\right| 
\preceq |\tau|^{-N}
$$
\textit{for  any positive} $N$.

We note that there is an integer $N_{0}$ independent of $\tau$ such that the total number of $k$ in the set $E_{3, k}$ is no more than $N_{0}$ for any fixed $\tau$. Therefore, it is easy to see that Lemma 3.1 follows from the following lemma (It can be found in Lemma 5.1 in \cite{bib6}). 

\textbf{Lemma 3.2.} \cite{bib6} \textit{Let} $k\in\mathbb{N}$, $0<\alpha<1$ \textit{and} $\beta>0$. \textit{For any nonnegative integers} $N$ \textit{and} $L$, \textit{we have that} 
$$
\left|\Gamma_{1, k}(\tau)\left(\frac{\partial}{\partial \tau}\right)^L \widehat{\mu_{\alpha, \beta, k}}(\tau)\right| \preceq 2^{k(L+1-\beta)} \Gamma_{1, k}(\tau) \min \left\{1,2^{-N k \alpha}\right\}, $$
$$\left|\Gamma_{2, k}(\tau)\left(\frac{\partial}{\partial \tau}\right)^L \widehat{\mu_{\alpha, \beta, k}}(\tau)\right| \preceq 2^{k(L+1-\beta)} \Gamma_{2, k}(\tau) \min \left\{1,\left|2^k \tau\right|^{-N}\right\}, $$
\begin{equation*}
	\begin{aligned}
		\left|\Gamma_{3, k}(\tau)\left(\frac{\partial}{\partial \tau}\right)^L \widehat{\mu_{\alpha, \beta, k}}(\tau)\right| &\preceq 2^{k(L+1-\beta)} \Gamma_{3, k}(\tau) \min \left\{1,|\tau|^{-1 / 2} 2^{-k / 2}\right\} \\
		&\approx 2^{k(L+1-\beta-\alpha / 2)} \Gamma_{3, k}(\tau). 
	\end{aligned}
\end{equation*}

By the same argument of the proof of Lemma 1 in \cite{bib6}, we have no difficulty to obtain the following result.

\textbf{Lemma 3.3.} \textit{Let} $k\in\mathbb{N}$, $0<\alpha<1$ \textit{and} $\beta>0$. \textit{Assume} $u>0$ \textit{and} $\frac{|x-y|(1+u)}{t}\le200$. \textit{For any nonnegative integer} $L$ \textit{and any multi-index} $\gamma$, \textit{we have that} 

$$
\begin{aligned}
	& \left|\frac{\partial^\gamma}{\partial y^\gamma} \frac{\partial^L}{\partial u^L} \widehat{\mu}_{\alpha, \beta}(\frac{|x-y|(1+u)}{t})\right| \\
	\preceq&\ t^{-|\gamma|}\left( \frac{|x-y|}{t}\right)  ^{L-\frac{L+|\gamma|}{1-\alpha}+\frac{\alpha-2+2 \beta}{2(1-\alpha)}}(1+u)^{|\gamma|+\frac{\alpha-2+2 \beta-2|\gamma|-2L}{2(1-\alpha)}} .
\end{aligned}
$$

\subsection{Derivative estimate of $\Omega_{\alpha, \beta, loc}(x, y, t)$}\label{subsec2}

The derivative estimate on $\Omega _{\alpha ,\beta ,loc}$ is crucial to the
proof of our main theorem. According to the result of Lemma 3.1, we have the
following Proposition 3.1, whose proof is similar to the proof of
Proposition 3 in a previous paper \cite{bib7}. However, the details of their proofs
seem slightly different, due to the appearance of parameter $t$ in
Proposition 3.1. For completeness, we shall give the details of proof.

\textbf{Proposition 3.1.} \textit{Let} $0<\alpha<1$, \textit{and} $\beta=\frac{n\alpha}{2}+\varepsilon$ \textit{for any fixed} $\varepsilon>0$. \textit{For any nonnegative integer} $L$ \textit{and any multi-index} $\gamma$, \textit{we have that}
$$ \left|\frac{\partial^\gamma}{\partial y^\gamma}\Omega_{\alpha, \beta, loc}(x, y, t)\right| \preceq t^{-n-|\gamma|}\left( \frac{|x-y|}{t}\right)  ^{-n+\frac{\varepsilon-|\gamma|}{1-\alpha}} \quad \textit{if}\quad \frac{|x-y|}{t}\le 1,
$$
and 
$$ \left|\frac{\partial^\gamma}{\partial y^\gamma}\Omega_{\alpha, \beta, loc}(x, y, t)\right| \preceq t^{-n-|\gamma|}\left( \frac{|x-y|}{t} \right) ^{-L} \quad \textit{if}\quad \frac{|x-y|}{t}> 1.
$$

\bigskip
\textbf{Proof.} By the Hadamard construction of a parametrix (see \cite{bib11,bib17}), we have that
$$\Omega_{\alpha, \beta, loc}(x, y, t)=\frac{2}{\pi t}\sum_{j=0}^{\infty} U_j(x, y)\int_{0}^{\infty} \tau \widehat {\mu}_{\alpha, \beta}(\tau/t)\frac{\left(\tau^2-|x-y|^2\right)_{+}^{j-\frac{n+1}{2}}}{\Gamma\left(j-\frac{n-1}{2}\right)}\psi(\tau)d\tau.$$ 
By the finite propagation speed property, we can easily see $\Omega_{\alpha, \beta, loc}(x, y, t)=0$ whenever $|x-y|\ge r_{0}$. In fact, since the support condition of $\psi$ and
$$\left(\tau^2-|x-y|^2\right)_{+}^{j-\frac{n+1}{2}}\neq0 $$
if $|x-y|\le|\tau|$, we must have
$$\Omega_{\alpha, \beta, loc}(x, y, t)=0$$
if $|x-y|\ge\frac{2r_{0}}{3}$.

We change variables $\tau=|x-y|(1+u)$ to obtain that
\begin{equation*}
	\begin{aligned}
		\Omega_{\alpha, \beta, loc}(x, y, t)
		=&\ t^{-1}\sum_{j=0}^{\infty}\left|x-y\right|^{-n+1+2j} U_j(x, y)V_{\alpha, \beta, j}(x, y, t)\\
		=&\sum_{j=0}^{\infty} t^{-n+2j} \left( \frac{|x-y|}{t}\right) ^{-n+1+2j} U_j(x, y)V_{\alpha, \beta, j}(x, y, t),
	\end{aligned}
\end{equation*}
where
\begin{equation*}
	\begin{aligned}
		V_{\alpha, \beta, j}(x, y, t)=&\frac{2}{\pi}\int_{0}^{\infty} \widehat {\mu}_{\alpha, \beta} (\frac{|x-y|(1+u)}{t}) (1+u)(2+u)^{j-\frac{n+1}{2}}\\
		&\times \psi(|x-y|(1+u))\frac{u_{+}^{j-\frac{n+1}{2}}}{\Gamma\left(j-\frac{n-1}{2}\right)}du,
	\end{aligned}
\end{equation*}
and for all multi-indices $\gamma$
$$\sum_{j=0}^{\infty}\left|\frac{\partial^\gamma}{\partial y^\gamma}U_j(x, y)\right|  <C$$
holds uniformly for $x, y\in\mathbb{M}$.
Now, we write
$$G_{j}(|x-y|, u, t)=F(u)  \widehat {\mu}_{\alpha, \beta} (\frac{|x-y|(1+u)}{t})\psi(|x-y|(1+u)),$$
where
$$F(u)=(1+u)(2+u)^{j-\frac{n+1}{2}}.$$ 
Then we rewrite $V_{\alpha, \beta, j}$ as 
$$V_{\alpha, \beta, j}(x, y, t)=\frac{2}{\pi}\int_{0}^{\infty}G_{j}(|x-y|, u, t)\frac{u_{+}^{j-\frac{n+1}{2}}}{\Gamma\left(j-\frac{n-1}{2}\right)}du.$$

By the Leibniz rule, for any multi-index $\gamma$ , we have that
\begin{equation*}
	\begin{aligned}
		& \frac{\partial^\gamma}{\partial y^\gamma}\Omega_{\alpha, \beta, loc}(x, y, t)\\
		=&\sum_{|K|\le|\gamma|}\sum_{j=0}^{\infty}C_{\gamma, K}\frac{\partial^K}{\partial y^K}(t^{-1}\left|x-y\right|^{-n+1+2j} U_j(x, y))\frac{\partial^{\gamma-K}}{\partial y^{\gamma-K}}V_{\alpha, \beta, j}(x, y, t),
	\end{aligned}
\end{equation*}
where $C_{\gamma, K}$ are constants depending on the multi-indices $K$ and $\gamma$. 

It is easy to see that 
$$
\begin{aligned}
	& \left|\frac{\partial^K}{\partial y^K}\left(t^{-1}|x-y|^{-n+1+2 j} U_j(x, y)\right) \frac{\partial^{\gamma-K}}{\partial y^{\gamma-K}} V_{\alpha, \beta, j}(x, y, t)\right| \\
	\preceq&\ t^{-1}|x-y|^{-n+1+2 j}\left|U_j(x, y)\right|\left|\frac{\partial^\gamma}{\partial y^\gamma} V_{\alpha, \beta, j}(x, y, t)\right|
\end{aligned}
$$ 
uniformly for multi-indices $K$ satisfying $|K|\le|\gamma|$ for each $j = 0, 1,\dots.$ 

So we only need to consider
the leading term $$t^{-1}\left|x-y\right|^{-n+1+2j} U_j(x, y)\frac{\partial^\gamma}{\partial y^\gamma}V_{\alpha, \beta, j}(x, y, t).$$
For $\frac {|x-y|}{t}> 1$, we can easily get it by Lemma 3.1. So we only need to consider the case of $\frac {|x-y|}{t}\le 1$ below.

Using the Leibniz rule and Lemma 3.3, it is easy to check that

\begin{eqnarray}
	&&\left|\frac{\partial^\gamma}{\partial y^\gamma} \frac{\partial^L}{\partial u^L} G_{j}(|x-y|, u, t)\right| \notag \\ 
	&\preceq& t^{-|\gamma|}\left( \frac{|x-y|}{t}\right)  ^{L-\frac{L+|\gamma|}{1-\alpha}+\frac{\alpha-2+2 \beta}{2(1-\alpha)}}(1+u)^{|\gamma|+\frac{\alpha-2+2 \beta-2|\gamma|-2L}{2(1-\alpha)}+j-\frac{n-1}{2} }.\label{4}
\end{eqnarray}

From the support condition of $\psi$, we have that
$$  G_{j}(|x-y|, u, t)=0 \quad \textup{if}\ u\ge \frac{2r_{0}}{3|x-y|},$$
which combining the estimate \eqref{4} with $L=0$ yields that
\begin{eqnarray}
	&&\left|\frac{\partial^\gamma}{\partial y^\gamma} V_{\alpha, \beta, j}(x, y, t)       \right| \notag \\ 
	&\preceq& t^{-|\gamma|}\left( \frac{|x-y|}{t}\right) ^{\frac{\alpha-2+2 \beta}{2(1-\alpha)}-\frac{|\gamma|}{1-\alpha}}\int_{0}^{ \frac{2r_{0}}{3|x-y|}}(1+u)^{|\gamma|+\frac{\alpha-2+2 \beta}{2(1-\alpha)}-\frac{|\gamma|}{1-\alpha}+2j-n}du\label{5}
\end{eqnarray} 
in the case $j\ge\frac{n}{2}$ and $n$ is even or in the case $j\ge\frac{n+1}{2}$ and $n$ is odd.

In the following we continue to estimate the integral in \eqref{5} by
considering three different cases. Also, it suffices to estimate the leading
term in each case.

\textbf{CASE 1.} $|\gamma|+\frac{\alpha-2+2 \beta}{2(1-\alpha)}-\frac{|\gamma|}{1-\alpha}+2j-n=-1.$ 

In this case, we have that
\begin{equation*} 
	\begin{aligned}
		&\left|\frac{\partial^\gamma}{\partial y^\gamma} V_{\alpha, \beta, j}(x, y, t)       \right| \notag \\ 
		\preceq&\  t^{-|\gamma|}\left( \frac{|x-y|}{t}\right) ^{\frac{\alpha-2+2 \beta}{2(1-\alpha)}-\frac{|\gamma|}{1-\alpha}}\log\left(1+\frac{1}{|x-y|} \right),
	\end{aligned} 
\end{equation*}  
which leads to
\begin{equation*} 
	\begin{aligned}
		&\left| t^{-1}\left|x-y\right|^{-n+1+2j} U_j(x, y)\frac{\partial^\gamma}{\partial y^\gamma}V_{\alpha, \beta, j}(x, y, t)\right|\\
		\preceq&\ t^{-n+2j-|\gamma|}\left( \frac{|x-y|}{t}\right)  ^{-|\gamma|}|U_{j}(x,y)|\log\left(1+\frac{1}{|x-y|} \right).
	\end{aligned} 
\end{equation*}  

\textbf{CASE 2.} $|\gamma|+\frac{\alpha-2+2 \beta}{2(1-\alpha)}-\frac{|\gamma|}{1-\alpha}+2j-n<-1.$ 

In this case, we have that
$$\left|\frac{\partial^\gamma}{\partial y^\gamma} V_{\alpha, \beta, j}(x, y, t)       \right| \preceq\  t^{-|\gamma|}\left( \frac{|x-y|}{t}\right) ^{\frac{\alpha-2+2 \beta}{2(1-\alpha)}-\frac{|\gamma|}{1-\alpha}}$$
and
\begin{equation*} 
	\begin{aligned}
		&\left| t^{-1}\left|x-y\right|^{-n+1+2j} U_j(x, y)\frac{\partial^\gamma}{\partial y^\gamma}V_{\alpha, \beta, j}(x, y, t)\right|\\
		\preceq&\ t^{-n+2j-|\gamma|}\left( \frac{|x-y|}{t}\right)  ^{-n+1+2j+\frac{\alpha-2+2 \beta}{2(1-\alpha)}-\frac{|\gamma|}{1-\alpha}}|U_{j}(x,y)|.
	\end{aligned} 
\end{equation*}  

\textbf{CASE 3.} $|\gamma|+\frac{\alpha-2+2 \beta}{2(1-\alpha)}-\frac{|\gamma|}{1-\alpha}+2j-n>-1.$ 

In this case, we have that
\begin{equation*} 
	\begin{aligned}
		&\left|\frac{\partial^\gamma}{\partial y^\gamma} V_{\alpha, \beta, j}(x, y, t)       \right| \\
		\preceq&\  t^{n-2j-1-2|\gamma|-\frac{\alpha-2+2 \beta}{2(1-\alpha)}+\frac{|\gamma|}{1-\alpha}   }\left( \frac{|x-y|}{t}\right)  ^{-2j+n-|\gamma|-1},
	\end{aligned} 
\end{equation*} 
which leads to
\begin{equation*} 
	\begin{aligned}
		&\left| t^{-1}\left|x-y\right|^{-n+1+2j} U_j(x, y)\frac{\partial^\gamma}{\partial y^\gamma}V_{\alpha, \beta, j}(x, y, t)\right|\\
		\preceq&\ t^{-2|\gamma|-1-\frac{\alpha-2+2 \beta}{2(1-\alpha)}+\frac{|\gamma|}{1-\alpha} }\left( \frac{|x-y|}{t}\right)  ^{-|\gamma|} |U_{j}(x,y)|.
	\end{aligned} 
\end{equation*}  

The estimates on the above three cases yield, for $\frac {|x-y|}{t}\le 1$,
\begin{eqnarray}
	&&\left| t^{-1}\sum_{j=n/2}^{\infty}\left|x-y\right|^{-n+1+2j} U_j(x, y)\frac{\partial^\gamma}{\partial y^\gamma}V_{\alpha, \beta, j}(x, y, t)\right| \notag \\ 
	&\preceq& t^{-|\gamma|}\left( \frac{|x-y|}{t}\right)  ^{\frac{\alpha-2+2 \beta}{2(1-\alpha)}-\frac{|\gamma|}{1-\alpha}+1}\log\left(1+\frac{1}{|x-y|} \right).\label{6}
\end{eqnarray}  

Now we start to estimate terms whose subscripts are less than $j\le\frac{n-1}{2}$. Here, we need to consider two cases according to $n$ is odd or $n$ is even.

\textbf{CASE A.} $n$ is odd and $j\le\frac{n-1}{2}$.

This is a relatively simple case. Combining  the definition given in \eqref{3} and Lemma 3.3, we can directly obtain that
\begin{equation*} 
	\begin{aligned}
		&\left|\frac{\partial^\gamma}{\partial y^\gamma} V_{\alpha, \beta, j}(x, y, t)       \right| \\
		=&\frac{1}{\pi}\left|(-1)^{\frac{n-1}{2}-j+1}\frac{1}{\Gamma(1)}\frac{\partial^\gamma}{\partial y^\gamma} \frac{\partial^{\frac{n-1}{2}-j}}{\partial u^{\frac{n-1}{2}-j}} G_{j}(|x-y|, u, t)\right|_{u=0}\\
		\preceq&\ t^{-|\gamma|}\left( \frac{|x-y|}{t}\right)  ^{\frac{n-1}{2}-j+\frac{\alpha-2+2 \beta}{2(1-\alpha)}-\frac{\frac{n-1}{2}-j+|\gamma|}{1-\alpha}},
	\end{aligned} 
\end{equation*} 
which yields to
\begin{equation*} 
	\begin{aligned}
		&\left| t^{-1}\left|x-y\right|^{-n+1+2j} U_j(x, y)\frac{\partial^\gamma}{\partial y^\gamma}V_{\alpha, \beta, j}(x, y, t)\right|\\
		\preceq&\ t^{-n+2j-|\gamma|}\left( \frac{|x-y|}{t}\right) ^{ j-\frac{n-1}{2}+\frac{\alpha-2+2 \beta}{2(1-\alpha)}-\frac{\frac{n-1}{2}-j+|\gamma|}{1-\alpha}}|U_{j}(x,y)|.
	\end{aligned} 
\end{equation*}  

Therefore, for odd $n$, by the above estimate and \eqref{6} we see that if $\frac {|x-y|}{t}\le 1$
\begin{equation*} 
	\begin{aligned}
		&\left| t^{-1}\sum_{j=0}^{\infty}\left|x-y\right|^{-n+1+2j} U_j(x, y)\frac{\partial^\gamma}{\partial y^\gamma}V_{\alpha, \beta, j}(x, y, t)\right|\\
		\preceq&\ t^{-n-|\gamma|}\left( \frac{|x-y|}{t}\right)  ^{-\frac{n-1}{2}+\frac{\alpha-2+2 \beta}{2(1-\alpha)}-\frac{\frac{n-1}{2}+|\gamma|}{1-\alpha}}|U_{j}(x,y)|\\
		&+\sum_{j=1}^{\frac{n-1}{2}}t^{-n+2j-|\gamma|}\left( \frac{|x-y|}{t}\right) ^{ j-\frac{n-1}{2}+\frac{\alpha-2+2 \beta}{2(1-\alpha)}-\frac{\frac{n-1}{2}-j+|\gamma|}{1-\alpha}}|U_{j}(x,y)|\\
		&+ t^{-|\gamma|}\left( \frac{|x-y|}{t}\right) ^{\frac{\alpha-2+2 \beta}{2(1-\alpha)}-\frac{|\gamma|}{1-\alpha}+1}\log\left(1+\frac{1}{|x-y|} \right)\\
		\preceq&\ t^{-n-|\gamma|}\left( \frac{|x-y|}{t}\right) ^{-\frac{n}{2}+\frac{\beta-n/2-|\gamma|}{1-\alpha}}\left(1+O\left(\left( \frac{|x-y|}{t}\right) ^{1+\frac{1}{1-\alpha}}\right)\right).
	\end{aligned} 
\end{equation*}   

\textbf{CASE B.} $n$ is even and $j\le\frac{n-2}{2}$. 

Assume that $\omega_{0}$ is a $C^{\infty}$ nonnegative function supported in $|\tau|\le2$ and satisfies
$$\omega_{0}(\tau)\equiv1 \quad \textup{if}\ |\tau|\le1.$$
Let
$$\omega_{\infty}(\tau)=1-\omega_{0}(\tau).$$ 
Then we can obtain that 
$$ \frac{\partial^\gamma}{\partial y^\gamma} G_{j}(|x-y|, u, t)=G_{j, 1}^{(\gamma)}(|x-y|, u, t)+G_{j, 2}^{(\gamma)}(|x-y|, u, t),$$ 
where  
$$G_{j, 1}^{(\gamma)}(|x-y|, u, t)=\left(\frac{\partial^\gamma}{\partial y^\gamma} G_{j}(|x-y|, u, t)\right)\omega_{0}\left(u\left( \frac{|x-y|}{t}\right) ^{\frac{\alpha}{\alpha-1}}\right)$$  
and  
$$G_{j, 2}^{(\gamma)}(|x-y|, u, t)=\left(\frac{\partial^\gamma}{\partial y^\gamma} G_{j}(|x-y|, u, t)\right)\omega_{\infty}\left(u\left( \frac{|x-y|}{t}\right) ^{\frac{\alpha}{\alpha-1}}\right).$$ 
Therefore, we have that
\begin{equation*} 
	\begin{aligned}
		\left|\frac{\partial^\gamma}{\partial y^\gamma} V_{\alpha, \beta, j}(x, y, t)       \right| 
		\preceq& \left| \int_{0}^{\infty} G_{j, 1}^{(\gamma)}(|x-y|, u, t)\frac{u_{+}^{j-\frac{n+1}{2}}}{\Gamma\left(j-\frac{n-1}{2}\right)}du\right| \\
		&+\left| \int_{0}^{\infty} G_{j, 2}^{(\gamma)}(|x-y|, u, t)\frac{u_{+}^{j-\frac{n+1}{2}}}{\Gamma\left(j-\frac{n-1}{2}\right)}du\right|.
	\end{aligned} 
\end{equation*}  
Next, we will discuss the two terms on the right-hand side of the above inequality respectively.

For $G_{j, 2}^{(\gamma)}(|x-y|, u, t)$, by the Leibniz rule, we easily obtain that
\begin{equation*} 
	\begin{aligned}  
		&\left|\frac{\partial^K}{\partial y^K}(\psi(|x-y|(1+u))F(u))\frac{\partial^{\gamma-K}}{\partial y^{\gamma-K}} \widehat {\mu}_{\alpha, \beta} (\frac{|x-y|(1+u)}{t})\right|\\
		\preceq&\left|\psi(|x-y|(1+u))F(u)\left(  \frac{\partial^{\gamma}}{\partial y^{\gamma}} \widehat {\mu}_{\alpha, \beta} (\frac{|x-y|(1+u)}{t})             \right)\right|
	\end{aligned} 
\end{equation*} 
for any multi-index $K$ satisfying $|K|\le|\gamma|$. 

Using  Lemma 3.3 with $L=0$, we have that 

\begin{align*}
	&\left| \int_{0}^{\infty} G_{j, 2}^{(\gamma)}(|x-y|, u, t)\frac{u_{+}^{j-\frac{n+1}{2}}}{\Gamma\left(j-\frac{n-1}{2}\right)}du\right| \\
	\preceq&\int_{\left( \frac{|x-y|}{t}\right) ^{\frac{\alpha}{1-\alpha}}   }         ^{\frac{2r_{0}}{3|x-y|}}\left|\psi(|x-y|(1+u))(1+u)(2+u)^{j-\frac{n+1}{2}}\right|\\
	&\times\left|  \frac{\partial^{\gamma}}{\partial y^{\gamma}}  \widehat {\mu}_{\alpha, \beta} (\frac{|x-y|(1+u)}{t})\right| u_{+}^{j-\frac{n+1}{2}}du\\
	\preceq& \int_{\left( \frac{|x-y|}{t}\right) ^{\frac{\alpha}{1-\alpha}}}^{\frac{2r_{0}}{3|x-y|}}\left|  \frac{\partial^{\gamma}}{\partial y^{\gamma}}  \widehat {\mu}_{\alpha, \beta} (\frac{|x-y|(1+u)}{t})(2+u)^{j-\frac{n-1}{2}}u^{j-\frac{n+1}{2}}\right| du\\
	\preceq&\ t^{-|\gamma|}\left( \frac{|x-y|}{t}\right) ^{\frac{\alpha-2+2 \beta}{2(1-\alpha)}-\frac{|\gamma|}{1-\alpha}}\int_{\left( \frac{|x-y|}{t}\right) ^{\frac{\alpha}{1-\alpha}}}^{\frac{2r_{0}}{3|x-y|}}(1+u)^{|\gamma|-\frac{\alpha-2+2 \beta}{2(1-\alpha)}-\frac{|\gamma|}{1-\alpha}}\\
	&\times (2+u)^{j-\frac{n-1}{2}}u^{j-\frac{n+1}{2}}du\\
	\preceq&\ t^{-|\gamma|}\left( \frac{|x-y|}{t}\right) ^{\frac{\alpha-2+2 \beta}{2(1-\alpha)}-\frac{|\gamma|}{1-\alpha}}\int_{\left( \frac{|x-y|}{t}\right) ^{\frac{\alpha}{1-\alpha}}}^{1}u^{j-\frac{n+1}{2}}du\\
	&+ t^{-|\gamma|}\left( \frac{|x-y|}{t}\right) ^{\frac{\alpha-2+2 \beta}{2(1-\alpha)}-\frac{|\gamma|}{1-\alpha}}\int_{1}^{\frac{2r_{0}}{3|x-y|}} (1+u)^{|\gamma|-\frac{\alpha-2+2 \beta}{2(1-\alpha)}-\frac{|\gamma|}{1-\alpha}}(2+u)^{2j-n}du\\
	=&\ I+II.
\end{align*}

It is easy to check that  
$$I\preceq t^{-|\gamma|}\left( \frac{|x-y|}{t}\right) ^{\frac{n}{2}-1+\frac{-|\gamma|+j\alpha+\beta-n/2}{1-\alpha}}.$$ 

For $II$, we need to discuss it based on two different cases of power. If 
$$|\gamma|-\frac{\alpha-2+2 \beta}{2(1-\alpha)}-\frac{|\gamma|}{1-\alpha}+2j-n\le-1,$$  
then we have that
$$II\preceq t^{-|\gamma|}\left( \frac{|x-y|}{t}\right) ^{\frac{\alpha-2+2 \beta}{2(1-\alpha)}-\frac{|\gamma|}{1-\alpha}}\log\left(1+\frac{1}{|x-y|} \right).$$ 
If
$$|\gamma|-\frac{\alpha-2+2 \beta}{2(1-\alpha)}-\frac{|\gamma|}{1-\alpha}+2j-n>-1,$$    
then 
$$II\preceq t^{-|\gamma|}\left( \frac{|x-y|}{t}\right) ^{-|\gamma|-2j+n-1}.$$ 
Combining estimates on $I$ and $II$, we thus know that  
\begin{equation}
	\left| \int_{0}^{\infty} G_{j, 2}^{(\gamma)}(|x-y|, u, t)\frac{u_{+}^{j-\frac{n+1}{2}}}{\Gamma\left(j-\frac{n-1}{2}\right)}du\right| \preceq t^{-|\gamma|}\left( \frac{|x-y|}{t}\right) ^{\frac{n}{2}-1+\frac{-|\gamma|+j\alpha+\beta-n/2}{1-\alpha}}. \label{7}
\end{equation}  

Now, by the definition given in \eqref{3} and the estimate in \eqref{4}, we obtain that 
\begin{equation*} 
	\begin{aligned}  
		&\left| \int_{0}^{\infty} G_{j, 1}^{(\gamma)}(|x-y|, u, t)\frac{u_{+}^{j-\frac{n+1}{2}}}{\Gamma\left(j-\frac{n}{2}\right)}du\right|\\  
		=& \left|(-1)^{\frac{n}{2}-j}\int_{0}^{\infty}
		\frac{\partial^\gamma}{\partial y^\gamma} \frac{\partial^{\frac{n}{2}-j}}{\partial u^{\frac{n}{2}-j}} \left(G_{j, 1}(|x-y|, u, t)\right) \frac{u_{+}^{-\frac{1}{2}}}{\Gamma\left(\frac{1}{2}\right)}du\right|\\ 
		\preceq&\ t^{-|\gamma|}\left( \frac{|x-y|}{t}\right) ^{\frac{n}{2}-j+\frac{\alpha-2+2 \beta}{2(1-\alpha)}-\frac{\frac{n}{2}-j+|\gamma|}{1-\alpha}}\int_{ 0 }^{2 \left( \frac{|x-y|}{t}\right) ^{\frac{\alpha}{1-\alpha}}   }u^{-\frac{1}{2}}du.
	\end{aligned} 
\end{equation*} 
This further yields that
\begin{equation}
	\left| \int_{0}^{\infty} G_{j, 1}^{(\gamma)}(|x-y|, u, t)\frac{u_{+}^{j-\frac{n+1}{2}}}{\Gamma\left(j-\frac{n-1}{2}\right)}du\right| \preceq t^{-|\gamma|}\left( \frac{|x-y|}{t}\right) ^{\frac{n}{2}-1+\frac{-|\gamma|+j\alpha+\beta-n/2}{1-\alpha}}. \label{8}
\end{equation}  

Combining \eqref{7} and \eqref{8}, we now obtain that, for even $n$ and $j\le\frac{n-2}{2}$ and any multi-index $\gamma$, 
\begin{equation*} 
	\begin{aligned}
		&\left| t^{-1}\left|x-y\right|^{-n+1+2j} U_j(x, y)\frac{\partial^\gamma}{\partial y^\gamma}V_{\alpha, \beta, j}(x, y, t)\right|\\
		\preceq&\ t^{-n+2j-|\gamma|} \left( \frac{|x-y|}{t}\right) ^{2j-\frac{n}{2}+\frac{-|\gamma|+j\alpha+\beta-n/2}{1-\alpha}}|U_{j}(x,y)|.
	\end{aligned} 
\end{equation*}    
Similar to Case A, in Case B we obtain that, for even $n$,
\begin{equation*} 
	\begin{aligned}
		&\left| t^{-1}\sum_{j=0}^{\infty}\left|x-y\right|^{-n+1+2j} U_j(x, y)\frac{\partial^\gamma}{\partial y^\gamma}V_{\alpha, \beta, j}(x, y, t)\right|\\ 
		\preceq&\ t^{-n-|\gamma|}\left( \frac{|x-y|}{t}\right) ^{-\frac{n}{2}+\frac{\beta-n/2-|\gamma|}{1-\alpha}}\left(1+O\left(\left( \frac{|x-y|}{t}\right) ^{1+\frac{1}{1-\alpha}}\right)\right).
	\end{aligned} 
\end{equation*} 

The proof of Proposition 3.1 is completed. Now, we have established the estimate on kernel of $T_{\alpha, \beta, loc}$. The estimate of $T^{*}_{\alpha, \beta, loc}(f)$ is included in the proof of Theorem 1.

\section{Estimate of $T^{*}_{\alpha, \beta, \infty}(f)$}\label{sec4}

In this section, we establish the following proposition.

\textbf{Proposition 4.1.} \textit{Let} $0 <\alpha< 1$ \textit{and} $\beta>0$, \textit{we have that for any} $0<p<1$
\begin{equation*}
	\left\Vert T^{*}_{\alpha, \beta, \infty}(f)\right\Vert_{L^{p}(\mathbb{M})}\preceq \left\Vert f\right\Vert_{H^{p}(\mathbb{M})}. 
\end{equation*}

To prove this proposition, we need the following lemma in \cite{bib5}.

\textbf{Lemma 4.1.} \cite{bib5} \textit{Suppose that f is a differentiable function on} $\left[0, \sigma\right]$. \textit{Then we have that}
$$
\sup _{0<t \leq \sigma}|f(t)| \preceq\left(b \int_0^\sigma t^{\varepsilon}\left|f^{\prime}(t)\right|^2 d t\right)^{1 / 2}+\left(\frac{1}{b} \int_0^\sigma|f(t)|^2 \frac{d t}{t^{\varepsilon}}\right)^{1 / 2}+|f(0)|
$$
for any positive number $b$ and real number $\varepsilon$.

\textbf{Proof of Proposition 4.1.} We recall that
$$\widehat{\mu}_{\alpha, \beta}(\tau)=\sum_{k=0}^{\infty}\widehat{\mu}_{\alpha, \beta, k}(\tau),$$
then we have that
$$ \Omega_{\alpha, \beta, \infty}(x, y, t)=\sum_{k=0}^{\infty}\Omega_{\alpha, \beta, k, \infty}(x, y, t), $$
where
$$ \Omega_{\alpha, \beta, k, \infty}(x, y, t)=\frac{2}{\pi t} \int_{0}^{\infty} \widehat {\mu}_{\alpha, \beta, k}(\tau/t)\varphi(\tau)\cos (\tau \sqrt{\mathcal{L}})(x, y)d\tau. $$
The corresponding operator is
$$T_{\alpha, \beta, k, \infty}(f)(x, t)=\int_{\mathbb{M}}\Omega_{\alpha, \beta, k, \infty}(x, y, t)f(y)dy,$$
and the maximal operator is defined as
$$T^{*}_{\alpha, \beta, k, \infty}(f)(x)=\sup_{0<t\le\sigma}|T_{\alpha, \beta, k, \infty}(f)(x, t)|. $$

It suffices to show that there is a positive number $d$ so that
$$\left\Vert T^{*}_{\alpha, \beta, k, \infty}(f)\right\Vert_{L^{p}(\mathbb{M})}\preceq2^{-dk}\left\Vert f\right\Vert_{H^{p}(\mathbb{M})}.$$ 
By the definition, we have that
$$\int_{\mathbb{M}}\Omega_{\alpha, \beta, k, \infty}(x, y, t)f(y)dy=\sum_{j}\eta_{k}(t, \lambda_{j})\left\langle f, e_{j}\right\rangle e_{j}(x),$$
where 
$$\eta_{k}(t, \lambda_{j})=\frac{2}{\pi t} \int_{0}^{\infty} \widehat {\mu}_{\alpha, \beta, k}(s/t)\varphi(s)\cos (s |\lambda_{j}|)ds. $$

When $\lambda_{j}=0$, it is easy to check that
$$|\eta_{k}(t, 0)|\preceq1.$$

For $\lambda_{j}\neq0$, we take integration by parts $\ell$ times to obtain that
$$\eta_{k}(t, \lambda_{j})\approx t^{-1}(|\lambda_{j}|)^{-\ell} \int_{0}^{\infty} \frac{\partial^{\ell}}{\partial{s^{\ell}}}\left(\widehat{\mu_{\alpha, \beta, k}}(s/t)\varphi(s)\right)\cos(s\lambda_{j}+\frac{\pi\ell}{2})ds,$$
where $\ell=\left[n(\frac{1}{p}-\frac{1}{2})\right]+1$. 

Using the Leibniz rule for $\frac{\partial^{\ell}}{\partial{s^{\ell}}}\left(\widehat{\mu_{\alpha, \beta, k}}(s/t)\varphi(s)\right)$, as before, we only need to consider the leading term 
$$\varphi(s)\frac{\partial^{\ell}}{\partial{s^{\ell}}}\left(\widehat{\mu_{\alpha, \beta, k}}(s/t)\right).$$
We recall that 
$$\widehat{\mu}_{\alpha, \beta, k}(s/t)=2\int_{0}^{\infty}\frac{e^{i \lambda^\alpha}}{ \lambda^\beta}\phi (\frac{ \lambda}{2^{k}})\cos(\frac{s\lambda}{t})d\lambda.$$
Then
\begin{equation*}
	\begin{aligned}
		\frac{\partial^{\ell}}{\partial{s^{\ell}}}\left(\widehat{\mu_{\alpha, \beta, k}}(s/t)\right) &\approx t^{-\ell}\int_{0}^{\infty}\frac{\lambda^{\ell}e^{i \lambda^\alpha}}{ \lambda^\beta}\phi (\frac{ \lambda}{2^{k}})\cos(\frac{s\lambda}{t}+\dfrac{\pi\ell}{2})d\lambda\\
		&=\frac{2^{k(\ell-\beta+1)}}{t^{\ell}}\int_{0}^{\infty}\frac{\lambda^{\ell}e^{i (2^{k}\lambda)^\alpha}}{ \lambda^\beta}\phi (\lambda)\cos(\frac{2^{k}s\lambda}{t}+\dfrac{\pi\ell}{2})d\lambda.
	\end{aligned}
\end{equation*}
For the last integral, we continue to integrate by parts $N$ times on the variable $ \lambda $ so that
$$\left|\frac{\partial^{\ell}}{\partial{s^{\ell}}}\left(\widehat{\mu_{\alpha, \beta, k}}(s/t)\right)\right|\preceq \dfrac{2^{k(\ell-\beta+1)}2^{-kN(1-\alpha)}}{t^{\ell}s^{N}}t^{N}.$$ 
Note $ 0<\alpha<1 $, by choosing a sufficiently large $ N $, we have that
$$\left|\frac{\partial^{\ell}}{\partial{s^{\ell}}}\left(\widehat{\mu_{\alpha, \beta, k}}(s/t)\right)\right|\preceq 2^{-dk}t s^{-N}$$
for any $d>0$, where $ 0<t<\sigma. $

Since $\varphi(s)\neq0$ only if $|s|>r_{0}/3$, we obtain that there is a $C>0$ independent of $t$ such that
$$|\eta_{k}(t, \lambda_{j})|\le C 2^{-dk}|\lambda_{j}|^{-\ell}.$$
Note that
\begin{equation*}
	\begin{aligned}
		\frac{\partial}{\partial t}\eta_{k}(t, \lambda_{j})=&-\frac{2}{t^{2}} \int_{0}^{\infty} \widehat {\mu}_{\alpha, \beta, k}(s/t)\varphi(s)\cos (s |\lambda_{j}|)ds\\
		&+\frac{2}{t}\int_{0}^{\infty}\left(  \frac{\partial}{\partial t}\widehat {\mu}_{\alpha, \beta, k}(s/t)\right)\varphi(s)\cos (s |\lambda_{j}|)ds\\
		=&I+II.
	\end{aligned}
\end{equation*}
For $I$, a similar method used in the above estimate shows 
$$\left|\frac{2}{t^{2}} \int_{0}^{\infty} \widehat {\mu}_{\alpha, \beta, k}(s/t)\varphi(s)\cos (s |\lambda_{j}|)ds\right|\preceq 2^{-dk}|\lambda_{j}|^{-\ell}.$$
For $II$, we integrate by parts $\ell$ times. Then
\begin{equation*}
	\begin{aligned}
		&\frac{2}{t}\int_{0}^{\infty}\left(  \frac{\partial}{\partial t}\widehat {\mu}_{\alpha, \beta, k}(s/t)\right)\varphi(s)\cos (s |\lambda_{j}|)ds\\
		\approx&\frac{1}{t}|\lambda_{j}|^{-\ell}\int_{0}^{\infty}\varphi(s)\left( \frac{\partial^{\ell}}{\partial{s^{\ell}}} (\frac{\partial}{\partial t}\widehat {\mu}_{\alpha, \beta, k}(s/t))\right)\cos (s |\lambda_{j}|+\frac{\pi\ell}{2})ds.
	\end{aligned}
\end{equation*}
By the definition of $\widehat{\mu}_{\alpha, \beta, k}(s/t)$, we have
$$\frac{\partial}{\partial t}\widehat{\mu}_{\alpha, \beta, k}(s/t)\approx t^{-2}\int_{0}^{\infty}\frac{\lambda e^{i \lambda^\alpha}}{ \lambda^\beta}\phi (\frac{ \lambda}{2^{k}})s\cos(\frac{s\lambda}{t}+\dfrac{\pi}{2})d\lambda.$$
Then
\begin{equation*}
	\begin{aligned}
		& \frac{\partial^{\ell}}{\partial{s^{\ell}}} (\frac{\partial}{\partial t}\widehat {\mu}_{\alpha, \beta, k}(s/t))\\
		=&\ t^{-2} \int_{0}^{\infty}\frac{\lambda e^{i \lambda^\alpha}}{ \lambda^\beta}\phi (\frac{ \lambda}{2^{k}})\frac{\partial^{\ell}}{\partial{s^{\ell}}}\left[s\cos(\frac{s\lambda}{t}+\dfrac{\pi}{2})\right] d\lambda\\
		=&(\ell-1)t^{-1-\ell} \int_{0}^{\infty}\frac{\lambda^{\ell} e^{i \lambda^\alpha}}{ \lambda^\beta}\phi (\frac{ \lambda}{2^{k}})s\cos(\frac{\lambda}{t}+\dfrac{\pi\ell}{2})d\lambda\\
		&+	t^{-2-\ell} \int_{0}^{\infty}\frac{\lambda^{\ell+1} e^{i \lambda^\alpha}}{ \lambda^\beta}\phi (\frac{ \lambda}{2^{k}})s\cos(\frac{s\lambda}{t}+\dfrac{\pi(\ell+1)}{2})d\lambda.
	\end{aligned}
\end{equation*}

For the first integral in the inequality above, we treat it as before:
$$\left|(\ell-1)t^{-1-\ell} \int_{0}^{\infty}\frac{\lambda^{\ell} e^{i \lambda^\alpha}}{ \lambda^\beta}\phi (\frac{ \lambda}{2^{k}})s\cos(\frac{\lambda}{t}+\dfrac{\pi\ell}{2})d\lambda\right|\preceq  2^{-dk}t s^{-N}.$$

For the second integral, we can get that
\begin{equation*}
	\begin{aligned}
		&\left|t^{-2-\ell} \int_{0}^{\infty}\frac{\lambda^{\ell+1} e^{i \lambda^\alpha}}{ \lambda^\beta}\phi (\frac{ \lambda}{2^{k}})s\cos(\frac{s\lambda}{t}+\dfrac{\pi(\ell+1)}{2})d\lambda\right|\\
		\preceq&\  \dfrac{2^{k(\ell-\beta+2)}2^{-kN(1-\alpha)}}{t^{\ell+2}s^{N-1}}t^{N}.
	\end{aligned}
\end{equation*}
A similar method gets
$$	\left|t^{-2-\ell} \int_{0}^{\infty}\frac{\lambda^{\ell+1} e^{i \lambda^\alpha}}{ \lambda^\beta}\phi (\frac{ \lambda}{2^{k}})s\cos(\frac{s\lambda}{t}+\dfrac{\pi(\ell+1)}{2})d\lambda\right|
\preceq 2^{-dk}t s^{-N+1}.$$
Therefore, we have that
$$\left|\frac{2}{t}\int_{0}^{\infty}\left(  \frac{\partial}{\partial t}\widehat {\mu}_{\alpha, \beta, k}(s/t)\right)\varphi(s)\cos (s |\lambda_{j}|)ds\right|\preceq 2^{-dk}|\lambda_{j}|^{-\ell}.$$

Combining the above estimates, we get that
$$\left|\frac{\partial}{\partial t}\eta_{k}(t, \lambda_{j})\right|\preceq 2^{-dk}|\lambda_{j}|^{-\ell}.$$
\\
By the compactness of $\mathbb{M}$ and Hölder's inequality, we obtain that
$$\left\Vert T^{*}_{\alpha, \beta, k, \infty}(f)\right\Vert_{L^{p}(\mathbb{M})}\preceq \left\Vert T^{*}_{\alpha, \beta, k, \infty}(f)\right\Vert_{L^{2}(\mathbb{M})}.$$
Using Lemma 4.1 with $\varepsilon=0$ and $b=1$, we have that
\begin{equation*}
	\begin{aligned}
		\left\Vert T^{*}_{\alpha, \beta, k, \infty}(f)\right\Vert_{L^{2}(\mathbb{M})}^2\preceq& \left\| \left( \int_{0}^{\sigma}|\frac{\partial}{\partial t} T_{\alpha, \beta, k, \infty}(f)(t, \cdot)|^{2}dt\right) ^{\frac{1}{2}}\right\|_{L^{2}(\mathbb{M})}^2\\
		&+ \left\| \left( \int_{0}^{\sigma}| T_{\alpha, \beta, k, \infty}(f)(t, \cdot)|^{2}dt\right) ^{\frac{1}{2}}\right\|_{L^{2}(\mathbb{M})}^2.\\
		&+\left\Vert  \lim_{t\rightarrow 0+}T_{\alpha, \beta, k, \infty}(f)(t, \cdot) \right\Vert _{L^{2}(\mathbb{M}) }^{{}}.
	\end{aligned}
\end{equation*}
By Fubini's theorem and the Plancherel formula, the first term in the above inequality
\begin{equation*}
	\begin{aligned}
		& \left\| \left( \int_{0}^{\sigma}|\frac{\partial}{\partial t} T_{\alpha, \beta, k, \infty}(f)(t, \cdot)|^{2}dt\right) ^{\frac{1}{2}}\right\|_{L^{2}(\mathbb{M})}^2\\
		\preceq&\ 2^{-2dk}\int_{0}^{\sigma} \sum_{j}| \frac{\partial}{\partial t}\eta_{k}(t, \lambda_{j})|^{2}|\left\langle f, e_{j}\right\rangle |^{2}dt\\
		\preceq&\ 2^{-2dk}\left\| \mathfrak{R}_{-\ell}(f)  \right\|_{L^{2}(\mathbb{M})}^2.
	\end{aligned}
\end{equation*}
A similar method shows that
$$\left\| \left( \int_{0}^{\sigma}| T_{\alpha, \beta, k, \infty}(f)(t, \cdot)|^{2}dt\right) ^{\frac{1}{2}}\right\|_{L^{2}(\mathbb{M})}^2\preceq 2^{-2dk}\left\| \mathfrak{R}_{-\ell}(f)  \right\|_{L^{2}(\mathbb{M})}^2$$
and
$$\left\Vert  \lim_{t\rightarrow 0+}T_{\alpha, \beta, k, \infty}(f)(t, \cdot) \right\Vert _{L^{2}(\mathbb{M}) }^{2}\preceq 2^{-2dk}\left\| \mathfrak{R}_{-\ell}(f)  \right\|_{L^{2}(\mathbb{M})}^2.$$

By the Sobolev imbedding theorem, we may obtain that
$$\left\Vert T^{*}_{\alpha, \beta, k, \infty}(f)\right\Vert_{L^{p}(\mathbb{M})}\preceq\ 2^{-dk}\left\| f\right\| _{H^p(\mathbb{M})}.$$

The proof of Proposition 4.1 is completed.

\section{Proof of Theorem 1}\label{sec5}

 We first prove sufficiency of (i) in Theorem 1 here. In the previous section we have obtained the estimate of $T^{*}_{\alpha, \beta, \infty}(f)$. Next, we only need to use Proposition 3.1 in Section 3 to estimate $T^{*}_{\alpha, \beta, loc}(f)$. We obtain the following two propositions (Proposition 5.1 and Proposition 5.2). 

\textbf{Proposition 5.1.} \textit{Let} $0<\alpha<1$. \textit{For any} 
$$p=\frac{n}{n+m}, m=1, 2,\dots,$$ 
\textit{we have that}
$$\left\| T^{*}_{\alpha, \beta, loc}(f)\right\|_{L^{p, \infty}(\mathbb{M})}\preceq \left\| f\right\| _{H^p(\mathbb{M})},$$
\textit{where} $\beta=n\alpha(\frac{1}{p}-\frac{1}{2})=\frac{n\alpha}{2}+m\alpha.$

\textbf{Proof.}
By a standard atomic decomposition, it suffices to show that
$$\left\| T^{*}_{\alpha, \beta, loc}(a)\right\|_{L^{p, \infty}(\mathbb{M})}\preceq 1$$
uniformly on any atom $a$ (see \cite{bib22}). It is obvious if $a$ is an exceptional atom. Hence we will only consider $p$-atoms $a$. We assume supp$(a)=B(x_{0}, r)$.

We write
\begin{equation*}
	\begin{aligned}
		&\left| \left\lbrace x\in \mathbb{M} :| T^{*}_{\alpha, \beta, loc}(a)(x)|>\lambda     \right\rbrace\right|\\
		\preceq&\left| \left\lbrace |x-x_{0}|\le 8r :| T^{*}_{\alpha, \beta, loc}(a)(x)>\lambda     \right\rbrace\right|\\
		&+\left| \left\lbrace |x-x_{0}|> 8r  :| T^{*}_{\alpha, \beta, loc}(a)(x)|>\lambda    \right\rbrace\right|.
	\end{aligned}
\end{equation*}
By the Chebyshev inequality, H\"{o}lder's inequality and the $L^2(\mathbb{M})$ boundness of the maximal function $T^{*}_{\alpha, \beta}(f)$, we have that  
$$\left| \left\lbrace |x-x_{0}|\le 8r :| T^{*}_{\alpha, \beta, loc}(a)(x)|>\lambda     \right\rbrace\right|\preceq \frac{1}{\lambda^{p} } \int_{|x-x_{0}|\le 8r} \left|T^{*}_{\alpha, \beta}(a)(x)\right|^{p}dx\preceq \frac{1}{\lambda^{p} } .$$
Using the cancellation condition of $a$ and Proposition 3.1, we have that
\begin{equation*}
	\begin{aligned}
		&\left|\int_{\mathbb{M}}\Omega_{\alpha, \beta, loc}(x, y, t)a(y)dy\right|\\
		\preceq&\int_{|y-x_{0}|<r}|a(y)|\sup_{\xi\in B(x_{0}, r)}\sum_{|\gamma|=m} \left|\frac{\partial^\gamma}{\partial y^\gamma}\Omega_{\alpha, \beta, loc}(x, \xi, t)\right| |y-x_{0}|^{m}dy\\
		\preceq&\sup_{\xi\in B(x_{0}, r)}|x-\xi|^{-n-m}\int_{|y-x_{0}|<r}|a(y)||y-x_{0}|^{m}dy\\
		\preceq&\sup_{\xi\in B(x_{0}, r)}|x-\xi|^{-n-m}.
	\end{aligned}
\end{equation*}
So we obtain that
\begin{equation*}
	\begin{aligned}
		&\left| \left\lbrace |x-x_{0}|> 8r  :| T^{*}_{\alpha, \beta, loc}(a)(x)|>\lambda    \right\rbrace\right|\\
		\preceq&\left| \left\lbrace |x-x_{0}|> 8r  :|x-x_{0}|^{-n-m}>\lambda    \right\rbrace\right|\\
		\preceq&\ \frac{1}{\lambda^{p}} .
	\end{aligned}
\end{equation*}

Combining the above estimates, we have completed the proof of this proposition. 

Next, we estimate the following proposition.

\textbf{Proposition 5.2.} \textit{Let} $0<\alpha<1$. \textit{For any} 
$$\frac{n}{n+1}<p<1,$$ 
\textit{we have that}
$$\left\| T^{*}_{\alpha, \beta, loc}(f)\right\|_{L^{p, \infty}(\mathbb{M})}\preceq \left\| f\right\| _{H^p(\mathbb{M})},$$
\textit{where} $\beta=n\alpha(\frac{1}{p}-\frac{1}{2})=\frac{n\alpha}{2}+n\alpha(\frac{1}{p}-1).$

\textbf{Proof.} The proof is similar to the above, but mainly estimate the following measure
$$\left| \left\lbrace |x-x_{0}|> 8r  :| T^{*}_{\alpha, \beta, loc}(a)(x)|>\lambda    \right\rbrace\right|.$$
Recall that 
$$ \left|\frac{\partial^\gamma}{\partial y^\gamma}\Omega_{\alpha, \beta, loc}(x, y, t)\right| \preceq t^{-n-|\gamma|}\left( \frac{|x-y|}{t}\right)  ^{-n+\frac{\varepsilon-|\gamma|}{1-\alpha}}, $$
if $\beta=\frac{n\alpha}{2}+\varepsilon.$
Take 
$$\theta=n(\frac{1}{p}-1).$$
Clearly $0<\theta<1$. By Proposition 3.1 with $\gamma=1$ and $\gamma=0$, we have that
\begin{equation*}
	\begin{aligned}
		&\left|\Omega_{\alpha, \beta, loc}(x, y, t)-\Omega_{\alpha, \beta, loc}(x, x_{0}, t)\right|^{\theta}\left|\Omega_{\alpha, \beta, loc}(x, y, t)-\Omega_{\alpha, \beta, loc}(x, x_{0}, t)\right|^{1-\theta}\\
		\preceq&\ |y-x_{0}|^{\theta}\ t^{-n-\theta}\sup_{\xi\in B(x_{0}, r)}\left( \frac{|x-\xi|}{t} \right) ^{-n+\frac{n\alpha(\frac{1}{p}-1)}{1-\alpha}}\left( \frac{|x-\xi|}{t} \right)^{-\frac{\theta}{1-\alpha}}\\
		\preceq&\ |y-x_{0}|^{\theta}\sup_{\xi\in B(x_{0}, r)}|x-\xi|^{-n-\theta}.
	\end{aligned}
\end{equation*}
Thus, we have that
\begin{equation*}
	\begin{aligned}
		&\left|\int_{\mathbb{M}}\Omega_{\alpha, \beta, loc}(x, y, t)a(y)dy\right|\\
		\preceq&\sup_{\xi\in B(x_{0}, r)}|x-\xi|^{-n-\theta}\int_{|y-x_{0}|<r}|a(y)||y-x_{0}|^{\theta}dy\\
		\preceq&\sup_{\xi\in B(x_{0}, r)}|x-\xi|^{-n-\theta}.
	\end{aligned}
\end{equation*}
It further yields that
\begin{equation*}
	\begin{aligned}
		&\left| \left\lbrace |x-x_{0}|> 8r  :| T^{*}_{\alpha, \beta, loc}(a)(x)|>\lambda    \right\rbrace\right|\\
		\preceq&\left| \left\lbrace |x-x_{0}|> 8r  :| x-x_{0}|^{-n/p}>\lambda    \right\rbrace\right|\\
		\preceq&\ \frac{1}{\lambda^{p}} .
	\end{aligned}
\end{equation*}
The proof of the proposition is completed. 

Similar to the proofs of Proposition 5.1 and Proposition 5.2, we obtain the $H^{p}-L^{p, \infty}$ boundedness of $T^{*}_{\alpha, \beta, loc}$ for any
$$\frac{n}{n+m+1}<p<\frac{n}{n+m}$$
if
\[
\beta\ge n\alpha(1/p-1/2).
\]
Now we claim that this is also a necessary condition. We prove this issue using an
argument of contradiction. 

In fact, we have proved that (\cite{bib7})
\[
\left\Vert T_{\alpha,\beta }\left( f\right) \right\Vert _{L^{p}\left(
\mathbb{M}\right) }\preceq \left\Vert f\right\Vert _{H^{p}\left( \mathbb{M}\right) }
\]%
if and only $\beta \geq n\alpha \left( 1/p-1/2\right) .\ $Now we assume that
there are $p_{0}, \beta _{0}$ and $\delta >0$ such that 
\[
\beta _{0}=n\alpha \left( 1/p_{0}-1/2\right) -\delta 
\]%
and 
\[
\left\Vert T_{\alpha,\beta _{0}}^{ \ast }\left( f\right) \right\Vert
_{L^{p_{0},\infty }\left( \mathbb{M}\right) }\preceq \left\Vert f\right\Vert
_{H^{p_{0}}\left( \mathbb{M}\right) }.
\]%
Let 
\[
1/p_{1}=1/p_{0}-\frac{\delta }{n\alpha }.
\]%
Clearly $1>p_{1}>p_{0}$ if $\delta $ is sufficiently small. We also obtain  
\[
\beta _{0}=n\alpha \left( 1/p_{1}-1/2\right) .
\]
 
 Thus by the sufficient condition, we have that 
\[
\left\Vert T_{\alpha,\beta _{0}}^{ \ast }\left( f\right) \right\Vert
_{L^{p_{1},\infty }\left( \mathbb{M}\right) }\preceq \left\Vert f\right\Vert
_{H^{p_{1}}\left( \mathbb{M}\right) }.
\]%
An interpolation yields 
\[
\left\Vert T_{\alpha,\beta _{0}}^{ \ast }\left( f\right) \right\Vert
_{L^{p}\left( \mathbb{M}\right) }\preceq \left\Vert f\right\Vert _{H^{p}\left(
\mathbb{M}\right) }
\]%
for all $p\in \left( p_{0},p_{1}\right) .$ We observe that 
\[
\beta _{0}=n\alpha \left( 1/p_{1}-1/2\right) <n\alpha \left( 1/p-1/2\right), 
\]%
which contradicts the known result in \cite{bib7} since
\[
\left\Vert T_{\alpha,\beta_{0} }\left( f\right) \right\Vert _{L^{p}\left(
\mathbb{M}\right) }\le \left\Vert T_{\alpha,\beta _{0}}^{ \ast }\left( f\right) \right\Vert
_{L^{p}\left( \mathbb{M}\right) }.
\]
Then we complete the proof of (i) in Theorem 1.

For the $H^{p}-L^{p}$ boundedness of $T^{*}_{\alpha, \beta, loc}$, we have the following 
proposition. 

\textbf{Proposition 5.3.} \textit{Let} $0<\alpha<1$. \textit{For any} 
$$0<p<1,$$ 
\textit{we have that}
$$\left\| T^{*}_{\alpha, \beta, loc}(f)\right\|_{L^{p}(\mathbb{M})}\preceq \left\| f\right\| _{H^p(\mathbb{M})},$$
\textit{whenever} $\beta> n\alpha(\frac{1}{p}-\frac{1}{2}).$

It is easy to obtain the proof of the above proposition by (i) in Theorem 1 and an interpolation (see \cite{bib23}), we leave it to the reader.

Now, the proof of Theorem 1 is completed.

\section{Proofs of Theorem 2 and Theorem 3}\label{sec6}

\subsection{Proof of Theorem 2}\label{subsec6}

To prove Theorem 2, we write
$$T_{t}(f)(x)=t^{-\frac{\beta}{\alpha}}\left(\sum_{k=1}^{N}c_{k}A_{kt}^{\alpha}-I \right)(f)(x)$$
and
$$T^* (f)(x)=\sup_{0<t\le \sigma^{\alpha}}|T_{t}(f)(x)|,$$
where $\left\lbrace c_{k}\right\rbrace_{k=1}^{N}$ is to be chosen, and $I$ is the identity operator.
It suffices to show that
$$\left\| T^* (f)\right\|_{L^{p,\infty} (\mathbb{M})}\preceq\left\| f\right\|_{H^{p}_\beta (\mathbb{M})} $$
for any $f\in C^{\infty}\cap H^{p}_\beta (\mathbb{M})$.

Let $\Psi=1-\Phi$, the definition of $\Phi$ here is the same as above. Then
\begin{equation*}
	\begin{aligned}
		&t^{-\frac{\beta}{\alpha}}\sum_{k=1}^{N}c_{k}A_{kt}^{\alpha}(f)(x)\\
		=&t^{-\frac{\beta}{\alpha}}\sum_{\lambda_{j}}\sum_{k=1}^{N}c_{k}e^{ikt|\lambda_{j}|^{\alpha}}\Phi(t^{1/\alpha}|\lambda_{j}|)\left\langle f, e_j\right\rangle e_j(x)\\
		&+t^{-\frac{\beta}{\alpha}}\sum_{\lambda_{j}}\sum_{k=1}^{N}c_{k}e^{ikt|\lambda_{j}|^{\alpha}}\Psi(t^{1/\alpha}|\lambda_{j}|)\left\langle f, e_j\right\rangle e_j(x)\\
		=&\sum_{\lambda_{j}}\sum_{k=1}^{N}c_{k}\dfrac{e^{ikt|\lambda_{j}|^{\alpha}}}{(t^{1/\alpha}|\lambda_{j}|)^{\beta}}\Phi(t^{1/\alpha}  |\lambda_{j}|)|\lambda_{j}|^{\beta}\left\langle f, e_j\right\rangle e_j(x)\\
		&+\sum_{\lambda_{j}}\sum_{k=1}^{N}c_{k}\dfrac{e^{ikt|\lambda_{j}|^{\alpha}}}{(t^{1/\alpha}  |\lambda_{j}|)^{\beta}}\Psi(t^{1/\alpha}  |\lambda_{j}|)|\lambda_{j}|^{\beta}\left\langle f, e_j\right\rangle e_j(x).
	\end{aligned}
\end{equation*}

Thus  
\begin{equation*}
	\begin{aligned}
		&t^{-\frac{\beta}{\alpha}}\left(\sum_{k=1}^{N}c_{k}A_{kt}^{\alpha}-I \right)(f)(x)\\
		=&\sum_{\lambda_j}\dfrac{\sum_{k=1}^{N}c_{k}e^{ikt|\lambda_{j}|^{\alpha}}-1 }{(  t^{1/\alpha}  |\lambda_{j}|)^{\beta}  }\Psi(t^{1/\alpha}  |\lambda_{j}|)|\lambda_{j}|^{\beta}\left\langle f, e_j\right\rangle e_j(x)\\
		&+\sum_{k=1}^{N}c_{k}\sum_{\lambda_{j}}e^{ikt|\lambda_{j}|^{\alpha}}\dfrac{\Phi(t^{1/\alpha}  |\lambda_{j}|)}{( t^{1/\alpha}  |\lambda_{j}|)^{\beta}}|\lambda_{j}|^{\beta}\left\langle f, e_j\right\rangle e_j(x)\\
		&-\sum_{\lambda_{j}}\dfrac{\Phi(t^{1/\alpha}  |\lambda_{j}|)}{( t^{1/\alpha}  |\lambda_{j}|)^{\beta}}|\lambda_{j}|^{\beta}\left\langle f, e_j\right\rangle e_j(x).
	\end{aligned}
\end{equation*}
Write
$$\sum_{\lambda_j}|\lambda_{j}|^{\beta}\left\langle f, e_j\right\rangle e_j(x)=\mathfrak{R}_{\beta}(f)(x)=g(x).$$
Now we have that
\begin{equation*}
	\begin{aligned}
		&\sup_{0<t\le \sigma^{\alpha}}\left|   t^{-\frac{\beta}{\alpha}}\left(\sum_{k=1}^{N}c_{k}A_{kt}^{\alpha}(f)(x)-(f)(x)\right)  \right| \\
		=&\Omega_{1}(g)(x)+\sum_{k=1}^{N}\Omega_{2, k}(g)(x)+\Omega_{3}(g)(x),
	\end{aligned}
\end{equation*}
where 
$$\Omega_{1}(g)(x)=\sup_{0<s\le\sigma}\left|  \sum_{\lambda_j}\dfrac{\sum_{k=1}^{N}c_{k}e^{ik|s\lambda_{j}|^{\alpha}}-1 }{(  |s\lambda_{j}|)^{\beta}  }\Psi( |s\lambda_{j}|)\left\langle g, e_j\right\rangle e_j(x)    \right|,$$
$$\Omega_{2, k}(g)(x)\le \sum_{k=1}^{N}|c_{k}|\sup_{0<s\le\sigma}\left|  \sum_{\lambda_{j}}e^{ik|s\lambda_{j}|^{\alpha}}\dfrac{\Phi(|s\lambda_{j}|)}{(  |s\lambda_{j}|)^{\beta}}\left\langle g, e_j\right\rangle e_j(x)      \right|,$$
$$\Omega_{3}(g)(x)=\sup_{0<s\le\sigma}\left|    \sum_{\lambda_{j}}\dfrac{\Phi(|s\lambda_{j}|)}{(   |s\lambda_{j}|)^{\beta}}\left\langle g, e_j\right\rangle e_j(x)   \right|.$$
To prove the theorem, since we have showed (Theorem 1)
$$\left\| \Omega_{2, k}(g)\right\|_{L^{p, \infty}(\mathbb{M})}\preceq \left\| g\right\| _{H^p(\mathbb{M})}$$
for each $k$, it remains to show that
$$\left\| \Omega_{1}(g)\right\|_{L^{p, \infty}(\mathbb{M})}\preceq \left\| g\right\| _{H^p(\mathbb{M})},$$
$$\left\| \Omega_{3}(g)\right\|_{L^{p, \infty}(\mathbb{M})}\preceq \left\| g\right\| _{H^p(\mathbb{M})}.$$
For $\Omega_{3}(g)$, we can find its $H^{p}(\mathbb{M})$-$L^{p}(\mathbb{M})$ boundedness in \cite{bib19} (The detail proof can see the proofs of Proposition 1.1 and Proposition 1.3 in \cite{bib19}).

Next, we only need to estimate $\Omega_{1}(g)$. We choose a sequence of constants $\{c_{k}\}_{k=1}^N$ satisfying
\begin{equation}
	\left(
	\begin{array}{cccc}
		1 & 1 & \cdots & 1\\
		1 & 2 & \cdots & N\\
		\vdots & \vdots & \cdots & \vdots\\
		1 & 2^{N-1} & \cdots & N^{N-1}
	\end{array}
	\right)
	\left(
	\begin{array}{c}
		c_{1}\\
		c_{2}\\
		\vdots\\
		c_{N}
	\end{array}
	\right)
	=
	\left(
	\begin{array}{c}
		1\\
		0\\
		\vdots\\
		0
	\end{array}
	\right),\notag
\end{equation} 
and use the Taylor expansions, for $k=1, 2, \dots, N,$
$$c_{k}e^{ik|s\lambda_{j}|^{\alpha}}=\sum_{j=0}^{\infty}c_{k}\frac{(ik|s\lambda_{j}|^{\alpha})^{j}}{j!}$$
to see that
\begin{equation*}
	\begin{aligned}
		&\sum_{k=1}^{N}c_{k}e^{ik|s\lambda_{j}|^{\alpha}}-1 \\
		=&\left(\sum_{k=1}^{N}c_{k}\right)-1+\sum_{j=1}^{N-1}\frac{(i|s\lambda_{j}|^{\alpha})^{j}}{j!}\sum_{k=1}^{N}c_{k}k^{j}+\sum_{k=1}^{N}c_{k}\sum_{j=N}^{\infty}\frac{(ik|s\lambda_{j}|^{\alpha})^{j}}{j!}\\
		=&E(|s\lambda_{j}|),
	\end{aligned}
\end{equation*}
where
$$E(|s\lambda_{j}|)=\sum_{k=1}^{N}c_{k}\sum_{j=N}^{\infty}\frac{(ik|s\lambda_{j}|^{\alpha})^{j}}{j!}.$$
Thus, we have that
\begin{equation*}
	\begin{aligned}
		&\sum_{\lambda_j}\dfrac{\sum_{k=1}^{N}c_{k}e^{ik|s\lambda_{j}|^{\alpha}}-1 }{(  |s\lambda_{j}|)^{\beta}  }\Psi( |s\lambda_{j}|)\left\langle g, e_j\right\rangle e_j(x) \\
		=&\sum_{\lambda_{j}}B(|s\lambda_{j}|)\Psi( |s\lambda_{j}|)\left\langle g, e_j\right\rangle e_j(x) ,\\
	\end{aligned}
\end{equation*}
where
$$B(\xi)=\frac{E(\xi)}{|\xi|^{\beta}}.$$
Now, it is easy to check that $B$ is a $C^\infty$ function off the origin and satisfies
$$B(\xi)=c|\xi|^{\gamma}+C(\xi), \gamma\ge0 $$
with a constant $c$ and the function $C$ satisfying
$$\left|   \left(\frac{\partial^{k}}{\partial \xi^{k}}C\right)(\xi)  \right|=O(|\xi|^{\gamma+1-k})\quad \textup{for}\ |\xi|\le2.$$
By the proofs of Proposition 1.1 and Proposition 1.2 in \cite{bib19}, we can obtain that
$$\left\| \sup_{s>0}\left|   \sum_{\lambda_{j}}B(|s\lambda_{j}|)\Psi( |s\lambda_{j}|)\left\langle g, e_j\right\rangle e_j(\cdot)   \right|    \right\|_{L^{p, \infty}(\mathbb{M})}\preceq \left\| g\right\| _{H^p(\mathbb{M})}.$$
The theorem is proved.

\subsection{Proof of Theorem 3}\label{subsec6}

By changing variables, we have
\begin{equation*}
	\begin{aligned}
		&kt^{-k}\int_{0}^{t}(t-s)^{k-1}e^{is|\lambda_{j}|^{\alpha}}ds\\
		=&k\int_{0}^{1}(1-r)^{k-1}e^{itr|\lambda_{j}|^{\alpha}}dr.
	\end{aligned}
\end{equation*}
Therefore
$$I_{k,\alpha}(|\mathcal{L}|)(f)(x,t)=k\sum_{\lambda_j}\int_{0}^{1}(1-r)^{k-1}e^{itr|\lambda_{j}|^{\alpha}}dr\left\langle f, e_j\right\rangle e_j(x).$$
Let $\Psi=1-\Phi$, the definition of $\Phi$ is as above. Then
$$I_{k,\alpha}(|\mathcal{L}|)(f)=I_{k,\alpha,loc}(|\mathcal{L}|)(f)+I_{k,\alpha,\infty}(|\mathcal{L}|)(f),$$
where
$$I_{k,\alpha,loc}(|\mathcal{L}|)(f)(x,t)=k\sum_{\lambda_j}\int_{0}^{1}(1-r)^{k-1}e^{itr|\lambda_{j}|^{\alpha}}dr \Psi(t|\lambda_{j}|^{\alpha})\left\langle f, e_j\right\rangle e_j(x),$$
$$I_{k,\alpha,\infty}(|\mathcal{L}|)(f)(x,t)=k\sum_{\lambda_j}\int_{0}^{1}(1-r)^{k-1}e^{itr|\lambda_{j}|^{\alpha}}dr \Phi(t|\lambda_{j}|^{\alpha})\left\langle f, e_j\right\rangle e_j(x).$$
For $I_{k,\alpha,loc}(|\mathcal{L}|)(f)$, by the Taylor expansions, we have that
$$
I_{k,\alpha,loc}(|\mathcal{L}|)(f)(x,t)=\sum_{\lambda_{j}}(1+G(t|\lambda_{j}|^{\alpha}))\Psi(t|\lambda_{j}|^{\alpha})\left\langle f, e_j\right\rangle e_j(x),
$$
where
$$G(t|\lambda_{j}|^{\alpha})=k\sum_{n=1}^{\infty}\int_{0}^{1}\frac{(1-r)^{k-1}r^{n}}{n!}dr(t|\lambda_{j}|^{\alpha})^{n}.$$
A similar method of estimating $\Omega_{1}(g)$ in Theorem 2 shows
$$\left\| \sup_{t>0}\left|   \sum_{\lambda_{j}}(1+G(t|\lambda_{j}|^{\alpha}))\Psi( t|\lambda_{j}|^{\alpha})\left\langle f, e_j\right\rangle e_j(\cdot)   \right|    \right\|_{L^{p,\infty}(\mathbb{M})}\preceq \left\| f\right\| _{H^p(\mathbb{M})}.$$
We have that
\begin{equation*}
	\begin{aligned}
		&k\int_{0}^{1}(1-r)^{k-1}e^{itr|\lambda_{j}|^{\alpha}}dr\\
		=&k\int_{0}^{1}(1-r)^{k-1}\cos(tr|\lambda_{j}|^{\alpha})dr\\
		+&it|\lambda_{j}|^{\alpha} \int_{0}^{1}(1-r)^{k}\cos(tr|\lambda_{j}|^{\alpha})dr\\
		\approx&\ V_{k-\frac{1}{2}}(t|\lambda_{j}|^{\alpha})+it|\lambda_{j}|^{\alpha} V_{k+\frac{1}{2}}(t|\lambda_{j}|^{\alpha}),
	\end{aligned}
\end{equation*}
where
$$V_{\nu}(z)=\frac{J_{\nu}(z)}{z^\nu}$$
and $J_{\nu}(z)$ is the Bessel function of order $\nu.$ Therefore,
 we notice that the leading term of $I_{k,\alpha,\infty}(|\mathcal{L}|)$ is
$$\sum_{\lambda_{j}} \frac{e^{it|\lambda_{j}|^{\alpha}}}{(t|\lambda_{j}|^{\alpha})^{k}} \Phi(t|\lambda_{j}|^{\alpha})\left\langle f, e_j\right\rangle e_j(x).$$
According to Theorem 1, we have completed the proof of Theorem 3.

\bigskip

\end{document}